\documentclass[11pt]{article}

\usepackage[a4paper,margin=1.05in]{geometry}
\usepackage{amsmath,amssymb,amsthm,mathtools}
\usepackage{enumitem}
\usepackage{hyperref}
\usepackage{mathrsfs}
\usepackage[T1]{fontenc}
\usepackage{lmodern}
\usepackage{microtype}

\hypersetup{colorlinks=true,linkcolor=blue,citecolor=blue,urlcolor=blue}

\newtheorem{theorem}{Theorem}[section]
\newtheorem{proposition}[theorem]{Proposition}
\newtheorem{lemma}[theorem]{Lemma}
\newtheorem{corollary}[theorem]{Corollary}
\newtheorem{mainthm}{Theorem}

\theoremstyle{definition}
\newtheorem{definition}[theorem]{Definition}
\newtheorem{example}[theorem]{Example}
\theoremstyle{remark}
\newtheorem{remark}[theorem]{Remark}

\numberwithin{equation}{section}

\title{Principal Eigenvalues for Transport-Nonlocal Operators on Compact Manifolds: Ergodic Averages and Quotient Reduction}
\author{
    Xing Liang\thanks{School of Mathematical Sciences, and Institute of Geometry and Physics, University of
Science and Technology of China, Hefei 230026, China, and Hefei National Laboratory, Hefei 230088, China. Email: \texttt{xliang@ustc.edu.cn}.
    X. Liang is supported by NSFC12531008, NSFC12331006 and XDB0900100.}
    \and
    Lei Zhang\thanks{School of Mathematics and Statistics, Shaanxi Normal
    University, Xi'an, Shaanxi 710119, China. Email:
    \texttt{zhanglei890512@gmail.com}. L. Zhang is supported by the
    NSFC12471168, NSFC12171119 and GK202602001.}
    \and
    Niansen Zhang\thanks{School of Mathematical Sciences, University of Science and Technology of China,
    Hefei, Anhui 230026, China. Email:
    \texttt{neutralsentinel@outlook.com}.}
}
\date{}

\begin{document}
    \maketitle

    \begin{abstract}
        We study principal eigenvalues for transport--nonlocal operators on
        connected compact Riemannian manifolds. Transport changes the natural
        spectral threshold from the pointwise maximum of the potential to a
        maximal invariant average. We establish a necessary and sufficient
        criterion based on a resolvent spectral radius and derive sufficient
        existence conditions for uniquely ergodic flows and, through quotient
        reduction, for flows arising from compact abelian Lie group actions.
        Examples on the two-dimensional torus show that, for weak nonlocal
        coupling, transport can create or suppress principal eigenvalues, while
        sphere rotations illustrate how the quotient criterion accommodates
        singular orbits.
    \end{abstract}



    \section{Introduction}

    This paper studies principal eigenvalues for transport--nonlocal operators on compact manifolds.  We work throughout on a connected compact \(C^\infty\) Riemannian manifold \((M,g)\) without boundary, and \(m\) denotes the normalized Riemannian volume.  Let \(V\) be a \(C^1\) vector field on \(M\), let \((\Phi_t)_{t\in\mathbb{R}}\) be its flow, and let \(a\in C(M;\mathbb{R})\).

    Given a kernel \(K\), define the nonlocal operator
    \[
    (\mathcal{K}u)(x):=\int_M K(x,y)u(y)\,\mathrm{d}m(y),
    \qquad x\in M.
    \]
    The kernel is assumed throughout to satisfy the following condition.
    \begin{description}[leftmargin=2.8em]
        \item[(K)] \(K\in C(M\times M;\mathbb{R})\), \(K\ge0\), and there exist constants \(\rho_0>0\) and \(\kappa_0>0\) such that
        \[
        K(x,y)\ge\kappa_0
        \qquad\text{whenever }d_g(x,y)<\rho_0.
        \]
    \end{description}
    This standing assumption will not be repeated in later statements.
    We consider
    \[
    \mathcal{L}:=V+a+\mathcal{K}
    \]
    on \(C(M;\mathbb{R})\), where \(a\) acts by multiplication and the domain is determined by the transport generator.  The precise domain, order notation, and complexification conventions are fixed in Section~\ref{sec:abstract-criterion}.

    \begin{definition}\label{def:principal-eigenvalue}
        The spectral bound of \(\mathcal{L}\) is
        \[
        s(\mathcal{L})
        :=\sup\{\operatorname{Re}z:z\in\sigma(\mathcal{L})\},
        \]
        where the spectrum is taken after complexification.  A real number \(\lambda\) is called a principal eigenvalue of \(\mathcal{L}\) if \(\lambda=s(\mathcal{L})\) and there exists \(\varphi\) in the domain of \(\mathcal{L}\) such that
        \[
        \mathcal{L}\varphi=\lambda\varphi,
        \qquad
        \min_M\varphi>0.
        \]
        In this case \((\lambda,\varphi)\) is called a principal eigenpair.
    \end{definition}

    The usual nonlocal dispersal operator
    \[
    \int_M K(x,y)(u(y)-u(x))\,\mathrm{d} m(y)+c(x)u(x)
    \]
    is included in this framework after the loss term is absorbed into the potential $a$.

    \paragraph{Background and motivation.}
    When they exist, principal eigenvalues characterize leading exponential growth rates and persistence thresholds in positive linear dynamics. In heterogeneous media, they are used to characterize KPP spreading speeds~\cite{FreidlinGartner1979SMD, Weinberger2002JMB, ShenZhang2010JDE}, establish the existence of pulsating fronts~\cite{BerestyckiHamelNadirashvili2005JEMS, CovilleDavilaMartinez2013AIHPANL, ShenZhang2012PAMS}, and derive persistence criteria for population dynamics~\cite{RawalShen2012JDDE}. For uniformly elliptic operators on compact domains, the maximum principle and Krein--Rutman theorem yield existence and simplicity results. On more general domains, this framework extends to generalized principal eigenvalues via maximum-principle characterizations~\cite{BerestyckiNirenbergVaradhan1994CPAM, BerestyckiRossi2006JEMS}.

    For nonlocal dispersal problems, the spectral bound may not be an eigenvalue. This can occur even for
    \[
    u\longmapsto a(x)u+\int_M K(x,y)u(y)\,\mathrm{d}m(y).
    \]
    For scalar nonlocal dispersal equations, Bates and Zhao~\cite{BatesZhao2007JMAA} developed a principal eigenvalue theory for bounded domains using resolvent-positive operator theory (see, e.g.,~\cite{Thieme1998DCDS}).
    Shen and Zhang~\cite{ShenZhang2010JDE} developed a principal eigenvalue framework for space-periodic nonlocal dispersal operators via the resolvent-perturbation method for positive semigroups introduced by B\"urger~\cite{Burger1988MZ}.
    Via generalized Krein--Rutman theory~\cite{EdmundsPotterStuart1972PRSA, Nussbaum1981LNM}, Coville~\cite{Coville2010JDE} obtained sufficient conditions for realizing the principal eigenvalue by a continuous positive eigenfunction.
    When such an eigenfunction does not exist, measure-valued formulations provide natural substitutes. Coville~\cite{Coville2013AML} proved the existence of positive measure-valued eigensolutions, including explicit singular examples.
    Li, Coville, and Wang~\cite{LiCovilleWang2017DCDS}, as well as Berestycki, Coville, and Vo~\cite{BerestyckiCovilleVo2016JFA}, subsequently developed the Collatz--Wielandt characterizations and the theory of generalized principal eigenvalues for nonlocal operators.

    These methods have also been extended to more general nonlocal operators and systems. B\"urger's positive-semigroup perturbation approach has been used to study principal spectrum points of time-periodic scalar equations~\cite{ShenXie2017DCDSB} and principal eigenvalues of cooperative systems~\cite{BaoShen2017PAMS}.
    Generalized Krein--Rutman arguments have been applied to partially degenerate nonlocal systems and equations with delay~\cite{Zhang2024ANS, LiangZhangZhao2019JDE}; related spectral problems for generalized convolution operators on the circle were investigated in~\cite{DingLiang2015SIAMJMA}.
    The resolvent-positive perturbation framework~\cite{Thieme1998DCDS} has also been used for age-structured models~\cite{KangRuan2022MA}, time-periodic cooperative systems~\cite{FengLiRuanXin2024SIAMJMA}, and delayed nonlocal dispersal systems~\cite{LiangZhangZhao2019JDE}.
    Collatz--Wielandt and related variational characterizations have been developed for cooperative systems with nonlocal or coupled dispersal~\cite{SuWangZhang2023JDE, SuLiLouWang2023JMPA}.
    More recently, generalized principal eigenvalues admitting such characterizations have been approximated by means of monotone upper and lower control systems~\cite{WangZhang2026JMPA}.
    Nonlocal diffusion operators without transport and their spectral
    properties have been studied on metric measure spaces and
    Riemannian manifolds
    \cite{RodriguezBernalSastreGomez2016PRSE,BandleGonzalezFontelosWolanski2018CPDE,RodriguezBernal2022NA}.

    The spectral behavior changes when nonlocal operators are coupled with transport, because the first-order generator is unbounded and long-time dynamics enter the spectral threshold.
    In one dimension, Coville and Hamel studied nonlocal operators with a drift of fixed sign~\cite{CovilleHamel2020NA}.
    More recently, Ducrot, Griette, and Liang treated periodic nonlocal problems with sign-changing advection, where resolvent-positive methods reveal uniqueness, multiplicity, and singular-measure alternatives~\cite{DucrotGrietteLiang2026Preprint}.

    The literature on transport--nonlocal operators remains limited, with results available mainly for one-dimensional or spatially periodic settings; the latter correspond to flat tori. To the best of our knowledge, the present paper is the first systematic study of principal eigenvalues for such operators on general compact Riemannian manifolds.

    \paragraph{Main results.}
    Our first result gives a general spectral-radius criterion for the
    existence of a principal eigenvalue on compact manifolds.  The criterion
    is then applied in two dynamical settings: uniquely ergodic flows, for
    which existence follows under either of two complementary hypotheses,
    and flows generated by dense one-parameter subgroups of compact abelian
    Lie groups acting on \(M\), for which a sufficient existence condition
    reduces to a scalar integral on the orbit quotient.  Before giving the
    precise statements, we introduce the transport growth bound and the
    resolvent formulation used in the spectral-radius criterion.

    Let \(\Phi_t\) be the flow generated by \(V\). The growth bound of
    the weighted transport semigroup associated with \(V+a\) is
    \begin{equation}\label{eq:omega}
        \omega(a,V)
        :=\lim_{t\to\infty}\frac1t\sup_{x\in M}\int_0^t a(\Phi_sx)\,\mathrm{d} s
        =\max_{\nu\in\mathcal M_\Phi}\int_M a\,\mathrm{d}\nu,
    \end{equation}
    where $\mathcal M_\Phi$ is the set of $\Phi_t$-invariant probability
    measures.  We refer to \(\omega(a,V)\) as the transport growth bound.
    It is determined by the maximal long-time averages of \(a\) along the
    flow.  When \(V\equiv0\), the corresponding quantity is the growth bound
    in the absence of transport, namely
    \(
    \omega(a,0)=\max_M a
    \).
    The characterization in \eqref{eq:omega} is a basic viewpoint in
    ergodic optimization for flows (see, e.g.,~\cite{MorroSantAnnaVarandas2020AHP}
    for maximizing measures in the continuous-time setting).  Thus, depending
    on the dynamics, this
    variational quantity may reduce to an average with respect to a unique
    invariant measure or to the maximum of averages along orbit closures;
    in general, it is the maximum over all invariant probability measures.

    For \(\lambda>\omega(a,V)\), let
    \[
    R_\lambda:=(\lambda-(V+a))^{-1}
    \]
    denote the resolvent of the transport--multiplication part.  A nonzero
    function \(u\) is an eigenfunction of \(\mathcal{L}\) with eigenvalue
    \(\lambda\) if and only if
    \[
    u=R_\lambda\mathcal{K}u.
    \]
    Thus the problem is governed by the spectral radius
    \(r(R_\lambda\mathcal{K})\), where \(r(T)\) denotes the spectral radius
    of a bounded operator \(T\).  The precise resolvent framework is
    developed in Section~\ref{sec:abstract-criterion}. The general criterion is as follows.

    \medskip
    \begin{mainthm}[Principal eigenvalues on compact manifolds]\label{thm:intro-A}
        The following conditions are equivalent:
        \begin{enumerate}[label=\textup{(\roman*)},leftmargin=2.6em]
            \item $\mathcal{L}$ has a principal eigenvalue;
            \item $s(\mathcal{L})>\omega(a,V)$;
            \item
            ${r}(R_{\lambda_0}\mathcal{K})>1$
            for some $\lambda_0>\omega(a,V)$.

        \end{enumerate}
        Whenever these conditions hold, there is a unique
        $\lambda_*>\omega(a,V)$ such that
        \[
        {r}(R_{\lambda_*}\mathcal{K})=1.
        \]
        Moreover, $\lambda_*=s(\mathcal{L})$ is an isolated, algebraically simple
        principal eigenvalue with a strongly positive eigenfunction, and no other
        spectral point lies in the half-plane
        $\{\operatorname{Re} z\ge\lambda_*\}$.
    \end{mainthm}
    The criterion established in Theorem~\ref{thm:intro-A}
    is rooted in classical positive perturbation theory and has also been
    used for nonlocal dispersal operators without transport
    (see, e.g.,~\cite{Thieme1998DCDS, BatesZhao2007JMAA, ShenZhang2010JDE}).
    In the present setting,
    the invariant dynamics enter through the transport growth bound
    \(\omega(a,V)\), while the nonlocal term enters through
    \(r(R_\lambda\mathcal K)\).
    Theorem~\ref{thm:spectral-radius} gives the spectral-radius construction,
    while Theorem~\ref{thm:strict-averaged-threshold} proves the full
    equivalence stated above.

    \medskip
    For the uniquely ergodic result below, let \(\mu\) denote the unique
    invariant probability measure of the \(V\)-flow \(\Phi_t\), and put
    \[
    \bar a_\mu:=\int_M a\,\mathrm{d}\mu.
    \]
    We use either of the following two conditions:
    \begin{description}[leftmargin=2.8em]
        \item[(BD)] the Birkhoff deviations of $a-\bar a_\mu$ along $\Phi_t$ satisfy a uniform one-sided lower bound, namely there exists $C_0>0$ such that, for all \(x\in M\) and \(t\ge0\),
        \[
        \int_0^t\bigl(a(\Phi_sx)-\bar a_\mu\bigr)\,\mathrm{d} s\ge -C_0;
        \]

        \item[(AC)] the invariant measure $\mu$ is absolutely continuous with respect to the Riemannian volume with bounded density, namely there exists $C_1>0$ such that
        \[
        \,\mathrm{d}\mu=\rho\,\mathrm{d} m,
        \qquad
        0\le \rho\le C_1\quad m\text{-a.e.}
        \]
    \end{description}

    \medskip
    \begin{mainthm}[Uniquely ergodic flows]\label{thm:intro-B}
        Suppose that the $V$-flow is uniquely ergodic with invariant probability measure $\mu$.  Then
        \[
        \omega(a,V)=\bar a_\mu=\int_M a\,\mathrm{d}\mu.
        \]
        If either \textup{(BD)} or \textup{(AC)} holds, then
        $\mathcal{L}=V+a+\mathcal{K}$ has a principal eigenvalue
        \[
        \lambda_*=s(\mathcal{L})>\bar a_\mu.
        \]
    \end{mainthm}

    Unique ergodicity identifies the transport growth bound. However, the
    nonlocal operator is defined using the Riemannian volume \(m\), rather than
    the invariant measure \(\mu\), and \(\mu\) may be singular with respect to
    \(m\). The theorem therefore uses two alternative hypotheses:
    \textup{(BD)} imposes a uniform one-sided lower bound on Birkhoff deviations
    along trajectories, whereas \textup{(AC)} requires \(\mu\) to have bounded
    density with respect to \(m\). Next, under the compact abelian group assumptions below, we treat flows that are uniquely ergodic on each orbit closure rather than on all of $M$.

    \medskip
    \begin{mainthm}[The case of compact abelian group actions]\label{thm:intro-C}
        Let $H$ be a compact abelian Lie group acting smoothly and isometrically on $M$.  Let
        \(
        \gamma:\mathbb{R}\to H
        \)
        be a smooth one-parameter subgroup with dense image, and define
        \[
        V(x)=\frac{\,\mathrm{d}}{\,\mathrm{d} t}\bigg|_{t=0}\gamma(t)\cdot x .
        \]
        The orbit closures, hereafter called leaves, are the compact \(H\)-orbits
        \[
        F_q=H\cdot x,
        \qquad q=\pi(x)\in Q:=M/H.
        \]
        For $q=\pi(x)$, let $\nu_q$ be the $H$-invariant probability measure on
        $F_q$ obtained by pushing forward normalized Haar measure on $H$ under
        the map $h\mapsto h\cdot x$.  Set
        \[
        \eta:=\pi_\# m,
        \qquad
        \bar a(q):=\int_{F_q}a\,\mathrm{d}\nu_q.
        \]
        Then $\bar a\in C(Q)$ and
        \[
        \omega(a,V)=\max_{q\in Q}\bar a(q).
        \]
        Moreover, if
        \begin{equation}\label{eq:quotient-integral-condition}
            \int_Q\frac{\,\mathrm{d}\eta(q)}{\lambda-\bar a(q)}\to+\infty
            \qquad \text{ as }\lambda\downarrow\max_Q\bar a,
        \end{equation}
        then $V+a+\mathcal{K}$ has a principal eigenvalue
        \[
        \lambda_*=s(V+a+\mathcal{K})>\max_{q\in Q}\bar a(q).
        \]
    \end{mainthm}
    The quotient $Q$ may fail to be a smooth manifold, and orbit closures of different dimensions are allowed.  The proof, including the Haar
    disintegration over these orbit closures, is given in
    Section~\ref{sec:group}. If the flow is generated by a one-parameter
    subgroup of a compact Lie group \(G\) acting smoothly and isometrically on
    \(M\), then the closure \(H\) of that subgroup is compact and abelian, and
    the result applies with this \(H\).

    \begin{remark}\label{rem:intro-zero-transport}
        The case without transport is included in Theorem~\ref{thm:intro-C} by taking the compact group action to be trivial.  Then every leaf is a single point, $Q=M$, $\eta=m$, and the leaf average is just $a$.  In particular, the quotient-integral condition reduces to
        \[
        \int_M\frac{\,\mathrm{d}m(x)}{\lambda-a(x)}\to+\infty
        \qquad(\lambda\downarrow\max_M a).
        \]
    \end{remark}

    \begin{remark}\label{rem:intro-time-change}
        Theorems~\ref{thm:intro-B} and~\ref{thm:intro-C} also admit positive
        time-change versions.  Let \(W=fV\), where \(f\in C^1(M)\) and
        \(f>0\).  The orbits are unchanged, and a \(V\)-invariant probability
        measure \(\mu\) is transformed into the corresponding \(W\)-invariant
        probability measure
        \[
        \,\mathrm{d}\mu_f
        =\frac{f^{-1}}{\int_M f^{-1}\,\mathrm{d}\mu}\,\mathrm{d}\mu.
        \]
        The same rule applies on each group orbit.  Precise statements are
        given in Corollaries~\ref{cor:time-changed-UE}
        and~\ref{cor:time-changed-group}.
    \end{remark}

    \paragraph{Applications and spectral effects.}
    The criterion established in Theorem~\ref{thm:intro-A} also yields a sharp threshold for the existence of a principal eigenvalue of $V+a+\varepsilon\mathcal K$ (Corollary~\ref{cor:critical-coupling}), which is then used in the creation and suppression results below.  On \(\mathbb{T}^2\), two
    complementary comparisons distinguish totally irrational transport,
    rational transport, and zero transport.

    \begin{enumerate}[label=\textup{(\roman*)},leftmargin=2.6em]
        \item \textbf{Totally irrational versus rational transport.}
        For a totally irrational direction, every orbit is dense in
        \(\mathbb{T}^2\), and the transport averages the potential over the
        whole torus. Thus, the transport growth bound is simply the global average
        of \(a\) (see Corollary~\ref{cor:totally-irrational} and
        \eqref{eq:T2-irrational}). By contrast, for a rational direction, the
        orbit closures are closed one-dimensional leaves: the potential is
        averaged only along these leaves, leaving a one-dimensional quotient
        variable (see Corollary~\ref{cor:partially-rational} and
        \eqref{eq:T2-rational-average}--\eqref{eq:T2-rational-integral-condition}).
        
        \item \textbf{Rational transport versus zero transport.}
        In the zero-transport case,
        \(\omega(a,0)=\max_{\mathbb{T}^2}a\), and the relevant
        integrability condition is governed by the behavior of
        \(\max_{\mathbb{T}^2}a - a\) near the maximizers of \(a\). Under rational transport,
        by contrast, both the growth bound and the corresponding quotient
        condition depend on the leaf-average function. For suitable
        potentials and sufficiently weak nonlocal coupling, this transition can
        either create or suppress a principal eigenvalue, as established in
        Lemma~\ref{lem:finite-no-advection} and Theorem~\ref{thm:advection-creates}
        (creation) and in Theorem~\ref{thm:advection-suppresses} (suppression).
    \end{enumerate}

    We also apply the quotient framework to an axial rotation on \(S^2\),
    where the regular orbit closures are latitude circles and the poles are
    zero-dimensional singular orbits.

    \paragraph{Organization of the paper.}
    Section~\ref{sec:abstract-criterion} proves the compact-manifold principal eigenvalue theorem and develops the spectral tools used later.  Section~\ref{sec:uniquely-ergodic} proves Theorem~\ref{thm:intro-B}, including its version for positive time changes.  Section~\ref{sec:group} proves Theorem~\ref{thm:intro-C} for compact abelian group actions, including the Haar disintegration over orbit closures.  Section~\ref{sec:applications} gives the torus and sphere applications, including the two-dimensional examples in which transport in a rational direction can create or suppress a principal eigenvalue.

    \section{Principal eigenvalues on compact manifolds}\label{sec:abstract-criterion}

    In this section, we establish the abstract spectral criterion used throughout the paper in the geometric setting of the introduction and under the standing kernel assumption \textup{(K)}. We begin by introducing some basic notation and definitions.

    Throughout this paper, \(C(M)\) denotes the real Banach lattice \(C(M;\mathbb{R})\), endowed with the supremum norm and its usual positive cone \(C(M)_+\).  The constant function with value one is denoted by \(\mathbf{1}\), and we write \(u\gg0\) when \(u\in C(M)_+\) and \(\min_Mu>0\).  For a real- or complex-valued function $h$, $\vert{}h\vert{}$ denotes its pointwise modulus, and order inequalities between real-valued functions are always assumed to hold pointwise. For standard background on Banach lattices and positive operators (see, e.g.,~\cite{Schaefer1974Book}).

    For complex spectral arguments, we pass to the complexification \(C(M;\mathbb{C})\).  Operators and their domains are extended complex-linearly and denoted by the same symbols.  Accordingly, \(\sigma(T)\), \(\sigma_p(T)\), \(\rho(T)\), \(r(T)\), and \(s(T)\), as well as algebraic multiplicities, refer to the complexification of \(T\).

    The vector field \(V\) generates the complete flow \(\Phi_t\).  For \(t\ge0\), the transport--multiplication semigroup on \(C(M)\) is
    \begin{equation}\label{eq:transport-semigroup}
        (T(t)u)(x)
        =\exp\left(\int_0^t a(\Phi_sx)\,\mathrm{d}s\right)u(\Phi_tx),
        \qquad u\in C(M),\quad x\in M.
    \end{equation}
    Its generator is \(\mathcal{A}=V+a\), with
    \[
    D(\mathcal{A})
    =\left\{u\in C(M):
    \lim_{t\to0}\frac{u\circ\Phi_t-u}{t}
    \text{ exists in }C(M)\right\},
    \qquad
    \mathcal{A}u=Vu+au,
    \]
    where \(Vu\) denotes the displayed limit.  The operator \(\mathcal{K}\) is the integral operator introduced in the introduction.  Since \(K\) is continuous on the compact space \(M\times M\), \(\mathcal{K}:C(M)\to C(M)\) is bounded, compact, and positive by the Arzel\`a--Ascoli theorem and assumption \textup{(K)}.  We set
    \[
    D(\mathcal{L})=D(\mathcal{A}),
    \qquad
    \mathcal{L}=\mathcal{A}+\mathcal{K}.
    \]

    The local lower bound in \textup{(K)} implies
    \begin{equation}\label{eq:K-one-lower}
        \mathcal{K}\mathbf{1}(x)
        =\int_M K(x,y)\,\mathrm{d}m(y)
        \ge \kappa_0\inf_{z\in M}m(B_g(z,\rho_0))=:\kappa_M>0,
        \qquad x\in M.
    \end{equation}
    Here \(B_g(x,\rho_0)=\{y\in M:d_g(x,y)<\rho_0\}\).

    For $z\in\mathbb{C}$ with $\operatorname{Re} z>\omega(a,V)$, let
    $R_z:=(z-\mathcal{A})^{-1}$ denote the resolvent of the complexification
    of $\mathcal{A}$ on $C(M;\mathbb{C})$.  Then, for every
    $h\in C(M;\mathbb{C})$,
    \begin{equation}\label{eq:resolvent-formula}
        (R_z h)(x)
        =\int_0^\infty
        \exp\left(-zt+\int_0^t a(\Phi_sx)\,\mathrm{d} s\right)
        h(\Phi_tx)\,\mathrm{d} t,
        \qquad x\in M.
    \end{equation}
    For real $\lambda>\omega(a,V)$, the operator $R_\lambda$ leaves
    $C(M)$ invariant and is a bounded positive operator on $C(M)$
    (see, e.g.,~\cite[Section~1.7]{Pazy1983Book}).  Throughout this section, set
    \begin{equation}\label{eq:P-lambda}
        P_\lambda:=R_\lambda\mathcal{K},
        \qquad \lambda>\omega(a,V).
    \end{equation}
    The operator family in \eqref{eq:P-lambda} yields the following fixed-point reduction.

    \begin{lemma}\label{lem:eigen-resolvent}
        Let $\lambda>\omega(a,V)$.
        Then
        \begin{equation}\label{eq:eigen-resolvent-equivalence}
            \lambda\in\sigma_p(\mathcal{A}+\mathcal{K})
            \quad\Longleftrightarrow\quad
            1\in\sigma_p(P_\lambda).
        \end{equation}
        More precisely, if $P_\lambda\varphi=\varphi$, then $(\mathcal{A}+\mathcal{K})\varphi=\lambda\varphi$. Conversely, if $(\mathcal{A}+\mathcal{K})\varphi=\lambda\varphi$, then $P_\lambda\varphi=\varphi$.
    \end{lemma}

    \begin{proof}
        If $P_\lambda\varphi=\varphi$, then $\varphi=R_\lambda\mathcal{K}\varphi$, hence
        \(
        (\lambda-\mathcal{A})\varphi=\mathcal{K}\varphi.
        \)
        Thus $(\mathcal{A}+\mathcal{K})\varphi=\lambda\varphi$.  Conversely, if $(\mathcal{A}+\mathcal{K})\varphi=\lambda\varphi$, then $(\lambda-\mathcal{A})\varphi=\mathcal{K}\varphi$, and applying $R_\lambda$ gives $\varphi=R_\lambda\mathcal{K}\varphi=P_\lambda\varphi$.  This proves \eqref{eq:eigen-resolvent-equivalence}.
    \end{proof}

    \begin{lemma}\label{lem:minorization}
        Fix $\Lambda>\omega(a,V)$.  There exist an integer $N\ge1$ and a constant $c_N>0$ such that, for every \(\lambda\in(\omega(a,V),\Lambda]\) and every \(h\in C(M)_+\setminus\{0\}\),
        \begin{equation}\label{eq:uniform-minorization}
            P_\lambda^N h(x)\ge c_N\int_M h\,\mathrm{d} m,
            \qquad
            \forall x\in M.
        \end{equation}
    \end{lemma}

    \begin{proof}
        Choose $0<\rho<\rho_0/2$.  By uniform continuity of the flow near $t=0$, there exists $\tau>0$ such that, for all \(x\in M\) and \(0\le t\le\tau\),
        \[
        d_g(\Phi_t x,x)<\rho .
        \]
        Take $\Lambda^+:=\max\{\Lambda,0\}$.  If $0\le t\le\tau$, $y\in B_g(x,\rho)$ and $\omega(a,V)<\lambda\le\Lambda$, then
        \[
        -\lambda t+\int_0^t a(\Phi_sx)\,\mathrm{d} s
        \ge -\bigl(\Lambda^++\|a\|_\infty\bigr)\tau,
        \qquad
        d_g(\Phi_tx,y)<2\rho<\rho_0.
        \]
        Hence the kernel assumption gives
        \(K(\Phi_t x,y)\ge\kappa_0\) for all \(x\in M\),
        \(0\le t\le\tau\), and \(y\in B_g(x,\rho)\).  Therefore, for every
        \(\lambda\in(\omega(a,V),\Lambda]\) and every \(h\in C(M)_+\),
        \[
        P_\lambda h(x)
        \ge k_\Lambda\int_{B_g(x,\rho)}h(y)\,\mathrm{d} m(y),
        \qquad
        \forall x\in M,
        \]
        where
        \[
        k_\Lambda:=\kappa_0\tau\exp\{- (\Lambda^++\|a\|_\infty)\tau\}>0.
        \]
        In other words, for every \(\lambda\in(\omega(a,V),\Lambda]\) and every
        \(h\in C(M)_+\),
        \[
        P_\lambda h\ge k_\Lambda A_\rho h,
        \]
        where
        \[
        (A_\rho h)(x):=\int_{B_g(x,\rho)}h(y)\,\mathrm{d} m(y),
        \qquad x\in M.
        \]
        The constant \(k_\Lambda\) is independent of \(\lambda\) and \(h\).

        Let $D=\operatorname{diam}M$.  Choose $N\ge2$ with $D/N<\rho/2$, and choose $0<\rho_1<\rho/4$.  It remains to show that the fixed power $A_\rho^N$ dominates the total mass.  Compactness gives
        \[
        v_{\rho_1}:=\inf_{z\in M}m(B_g(z,\rho_1))>0.
        \]

        We claim that, for every \(h\in C(M)_+\),
        \[
        A_\rho^N h(x)\ge v_{\rho_1}^{N-1}\int_M h\,\mathrm{d} m,
        \qquad
        \forall x\in M.
        \]
        By Tonelli's theorem (see, e.g.,~\cite[Theorem~2.37]{Folland1999Book}),
        \[
        A_\rho^N h(x)
        =\int_M h(y)\!
        \int_{M^{N-1}}\prod_{j=1}^N
        \mathbf 1_{\{d_g(\xi_{j-1},\xi_j)<\rho\}}
        \,\mathrm{d} m(\xi_1)\cdots\,\mathrm{d} m(\xi_{N-1})\,\mathrm{d} m(y),
        \]
        where $\xi_0=x$ and $\xi_N=y$.  Since $M$ is compact and connected, any two points $x,y\in M$ can be joined by a minimizing geodesic.  Choose points
        $x=z_0,z_1,\ldots,z_N=y$ on it with
        \[
        d_g(z_{j-1},z_j)\le d_g(x,y)/N\le D/N<\rho/2 .
        \]
        Restricting the inner variables to $\xi_j\in B_g(z_j,\rho_1)$, $1\le j\le N-1$, and using $\rho_1<\rho/4$, we obtain
        \[
        A_\rho^N h(x)\ge v_{\rho_1}^{N-1}\int_M h\,\mathrm{d} m,
        \]
        which proves the claimed estimate for $A_\rho^N$.

        Finally, the estimate $P_\lambda\ge k_\Lambda A_\rho$ can be iterated because both operators are positive.  Combining it with the preceding bound for $A_\rho^N$ gives
        \[
        P_\lambda^N h
        \ge k_\Lambda^N A_\rho^N h
        \ge k_\Lambda^N v_{\rho_1}^{N-1}\int_M h\,\mathrm{d} m.
        \]
        Hence \eqref{eq:uniform-minorization} holds with
        $c_N=k_\Lambda^N v_{\rho_1}^{N-1}$.
    \end{proof}

    \begin{lemma}\label{lem:monotone-spectral-radius}
        The map
        \(
        \lambda\longmapsto {r}(P_\lambda)
        \)
        is continuous and strictly decreasing on \((\omega(a,V),\infty)\).  Moreover
        \[
        {r}(P_\lambda)\longrightarrow0
        \qquad\text{as }\lambda\to+\infty.
        \]
    \end{lemma}

    \begin{proof}
        Since \(\mathcal{K}\) is compact and \(R_\lambda\) is bounded for
        \(\lambda>\omega(a,V)\), the operator
        \(P_\lambda\) is compact and positive.  We first show that
        \({r}(P_\lambda)>0\).  Indeed, Lemma~\ref{lem:minorization}, applied
        with any \(\Lambda>\lambda\), gives \(P_\lambda^N\mathbf{1}\ge c_N\mathbf{1}\).
        By positivity, for every \(n\ge1\),
        \[
        (P_\lambda^N)^n\mathbf{1}\ge c_N^n\mathbf{1}.
        \]
        Hence \(\|(P_\lambda^N)^n\|\ge c_N^n\).  Gelfand's formula (see, e.g.,~\cite[Chapter~I, \S4.6]{Kato1995Book}) yields that
        \[
        {r}(P_\lambda)^N
        ={r}(P_\lambda^N)
        =\lim_{n\to\infty}\|(P_\lambda^N)^n\|^{1/n}
        \ge c_N>0 .
        \]

        Since \(\lambda\mapsto P_\lambda\) is norm-continuous and \(P_\lambda\)
        is compact, the continuity of \(\lambda\mapsto{r}(P_\lambda)\) follows
        from standard perturbation theory for isolated spectral sets (see, e.g.,
        \cite[Chapter~IV, \S3.5]{Kato1995Book}).  The non-increasing monotonicity also follows from
        \cite[Theorem~2.5]{Thieme1998DCDS}; here we give a short direct argument
        which in fact yields strict monotonicity.  Let
        \(\omega(a,V)<\lambda_1<\lambda_2\).  The resolvent identity gives
        \[
        R_{\lambda_1}-R_{\lambda_2}
        =(\lambda_2-\lambda_1)R_{\lambda_1}R_{\lambda_2},
        \]
        and hence
        \[
        P_{\lambda_1}-P_{\lambda_2}
        =(\lambda_2-\lambda_1)R_{\lambda_1}R_{\lambda_2}\mathcal{K}\ge0 .
        \]
        By Gelfand's formula again, we have
        \[
        {r}(P_{\lambda_2})\le {r}(P_{\lambda_1}).
        \]

        By the Krein--Rutman theorem (see, e.g.,~\cite[Theorem~19.2]{Deimling1985Book}), choose
        \(0\ne\phi_2\in C(M)_+\) such that
        \(P_{\lambda_2}\phi_2={r}(P_{\lambda_2})\phi_2\).  Applying
        Lemma~\ref{lem:minorization} with \(\Lambda>\lambda_2\), we obtain
        \[
        {r}(P_{\lambda_2})^N\phi_2
        =P_{\lambda_2}^N\phi_2
        \ge c_N\Bigl(\int_M\phi_2\,\mathrm{d} m\Bigr)\mathbf{1},
        \]
        and therefore \(\phi_2\gg0\).  Since
        \(R_{\lambda_2}\mathcal{K}\phi_2=P_{\lambda_2}\phi_2
        ={r}(P_{\lambda_2})\phi_2\gg0\), the explicit formula for
        \(R_{\lambda_1}\) yields
        \(
        w:=(P_{\lambda_1}-P_{\lambda_2})\phi_2
        =(\lambda_2-\lambda_1)R_{\lambda_1}R_{\lambda_2}\mathcal{K}\phi_2\gg0.
        \)
        Thus
        \[
        P_{\lambda_1}\phi_2={r}(P_{\lambda_2})\phi_2+w .
        \]

        Again by the Krein--Rutman theorem, choose a nonzero positive
        eigenfunctional \(\psi_1^*\) of \(P_{\lambda_1}^*\) associated with
        \({r}(P_{\lambda_1})\).  We claim that \(\psi_1^*\) is strictly positive in the sense that
        \[
        \langle \psi_1^*,h\rangle>0
        \quad\text{for every } h\in C(M)_+\setminus\{0\}.
        \] Indeed, for \(0\ne h\in C(M)_+\),
        Lemma~\ref{lem:minorization}, applied with \(\Lambda>\lambda_1\), gives
        \(P_{\lambda_1}^N h\ge c_N(\int_M h\,\mathrm{d} m)\mathbf{1}\).  Hence
        \[
        \begin{aligned}
            \langle\psi_1^*,h\rangle
            &={r}(P_{\lambda_1})^{-N}
            \langle\psi_1^*,P_{\lambda_1}^N h\rangle  \\
            &\ge
            {r}(P_{\lambda_1})^{-N}c_N
            \Bigl(\int_M h\,\mathrm{d} m\Bigr)\langle\psi_1^*,\mathbf{1}\rangle>0 .
        \end{aligned}
        \]
        Pairing \(P_{\lambda_1}\phi_2={r}(P_{\lambda_2})\phi_2+w\) with
        \(\psi_1^*\), we get
        \[
        {r}(P_{\lambda_1})\langle\psi_1^*,\phi_2\rangle
        =
        {r}(P_{\lambda_2})\langle\psi_1^*,\phi_2\rangle
        +\langle\psi_1^*,w\rangle .
        \]
        It then follows from \(\langle\psi_1^*,\phi_2\rangle>0\) and
        \(\langle\psi_1^*,w\rangle>0\) that
        \({r}(P_{\lambda_1})>{r}(P_{\lambda_2})\).

        Finally, let \(a_{\max}=\max_M a\).  For \(\lambda>a_{\max}\), the resolvent formula gives
        \[
        \|R_\lambda\|
        \le\int_0^\infty e^{-(\lambda-a_{\max})t}\,\mathrm{d} t
        =\frac1{\lambda-a_{\max}}.
        \]
        Consequently,
        \(
        0\le {r}(P_\lambda)\le\|P_\lambda\|
        \le\|R_\lambda\|\,\|\mathcal{K}\|\longrightarrow0
        \ \text{as }\lambda\to+\infty.
        \)
    \end{proof}

    The next result is the main result of this section.

    \begin{theorem}\label{thm:spectral-radius}
        Suppose that there exists \(\lambda_0>\omega(a,V)\) such that
        \begin{equation}\label{eq:spectral-radius-condition}
            {r}(P_{\lambda_0})>1.
        \end{equation}
        Then the following assertions hold.
        \begin{enumerate}[label=\textup{(\roman*)},leftmargin=2.6em]
            \item There exists a unique \(\lambda_*>\omega(a,V)\) such that
            \begin{equation}\label{eq:critical-spectral-radius}
                {r}(P_{\lambda_*})=1.
            \end{equation}
            \item The spectral bound of \(\mathcal{L}=\mathcal{A}+\mathcal{K}\) is \(\lambda_*\).  More precisely,
            \(s(\mathcal{L})=\lambda_*>\omega(a,V)\), and there exists \(\varphi\in D(\mathcal{A})\), \(\varphi\gg0\), such that
            \(\mathcal{L}\varphi=\lambda_*\varphi\).
            \item The eigenvalue \(\lambda_*\) is algebraically simple and isolated.
            \item The eigenvalue \(\lambda_*\) is strictly dominant:
            \begin{equation}\label{eq:rightmost-spectrum}
                \sigma(\mathcal{L})\cap\{z\in\mathbb{C}:\operatorname{Re} z\ge\lambda_*\}=\{\lambda_*\}.
            \end{equation}
        \end{enumerate}
    \end{theorem}

    \begin{proof}

        \textup{(i)}
        By Lemma~\ref{lem:monotone-spectral-radius}, the function
        \(\lambda\mapsto{r}(P_\lambda)\) is continuous, strictly decreasing, and tends to zero as \(\lambda\to+\infty\).  Since \({r}(P_{\lambda_0})>1\), we may choose \(\lambda_1>\lambda_0\) such that \({r}(P_{\lambda_1})<1\).  The intermediate value theorem gives some \(\lambda_*\in(\lambda_0,\lambda_1)\) satisfying \eqref{eq:critical-spectral-radius}, and strict monotonicity gives uniqueness.  This proves \textup{(i)}.

        \medskip
        \textup{(ii)}
        Since \(P_{\lambda_*}\) is compact and positive, the Krein--Rutman theorem
        (see, e.g.,~\cite[Theorem~19.2]{Deimling1985Book}) yields a nonzero function
        \(\varphi\in C(M)_+\) and a nonzero positive functional
        \(\varphi^*\in C(M)^*_+\) such that
        \[
        P_{\lambda_*}\varphi=\varphi,
        \qquad
        P_{\lambda_*}^*\varphi^*=\varphi^*.
        \]
        In view of
        \(0\ne\varphi\in C(M)_+\), we have \(\int_M\varphi\,\mathrm{d} m>0\). Applying
        Lemma~\ref{lem:minorization} to \(\varphi\), and using
        \(P_{\lambda_*}^N\varphi=\varphi\), we obtain
        \(\varphi\ge c_N(\int_M\varphi\,\mathrm{d} m)\mathbf{1}\). Hence \(\varphi\gg0\).  Moreover,
        Lemma~\ref{lem:eigen-resolvent} gives \(\mathcal{L}\varphi=\lambda_*\varphi\).
        Thus \(\lambda_*\) is an eigenvalue of \(\mathcal{L}\) with a strongly positive
        eigenfunction.

        It remains to show that no spectral point lies to the right of \(\lambda_*\).
        Let \(z\in\mathbb{C}\) with \(\operatorname{Re} z>\omega(a,V)\).  By
        \eqref{eq:resolvent-formula} and the positivity of \(T(t)\),
        \[
        |R_zg|\le R_{\operatorname{Re} z}|g|,
        \qquad g\in C(M;\mathbb{C}).
        \]
        Since \(|\mathcal{K} h|\le\mathcal{K}|h|\), iteration and Gelfand's formula
        (see, e.g.,~\cite[Chapter~I, \S4.6]{Kato1995Book}) give the two useful
        domination estimates
        \begin{equation}\label{eq:complex-resolvent-domination}
            |R_z\mathcal{K} h|\le P_{\operatorname{Re} z}|h|,
            \qquad
            {r}(R_z\mathcal{K})\le{r}(P_{\operatorname{Re} z}),
            \qquad h\in C(M;\mathbb{C}).
        \end{equation}
        On \(D(\mathcal{A})\), the resolvent factorization is
        \begin{equation}\label{eq:resolvent-factorization}
            z-\mathcal{L}=(z-\mathcal{A})(I-R_z\mathcal{K}).
        \end{equation}
        If \(\operatorname{Re} z>\lambda_*\), then \({r}(P_{\operatorname{Re} z})<1\) by \textup{(i)}, and hence \({r}(R_z\mathcal{K})<1\). It follows that \(I-R_z\mathcal{K}\) is invertible, so \eqref{eq:resolvent-factorization} gives \(z\in\rho(\mathcal{L})\).  Therefore \(\sigma(\mathcal{L})\) has no point with real part larger than \(\lambda_*\).  Since we have already shown \(\lambda_*\in\sigma_p(\mathcal{L})\), it follows that \(s(\mathcal{L})=\lambda_*\).  This proves \textup{(ii)}.

        \medskip
        \textup{(iii)}
        We first prove that \(\lambda_*>\omega(a,V)\) is an isolated spectral point of
        \(\mathcal{L}\).
        Choose \(\delta_0>0\) such that
        \(\lambda_*-\delta_0>\omega(a,V)\).  Then \(z-\mathcal{A}\) is invertible
        whenever \(|z-\lambda_*|<\delta_0\).  By
        \eqref{eq:resolvent-factorization}, in this neighbourhood
        \(z-\mathcal{L}\) is invertible if and only if
        \(I-R_z\mathcal{K}\) is invertible.

        Set \(F(z):=I-R_z\mathcal{K}\) on \(C(M;\mathbb{C})\).  The family
        \(R_z\mathcal{K}\) is compact and analytic for
        \(\operatorname{Re} z>\omega(a,V)\).  Moreover, by \textup{(ii)}, \(F(z)\) is
        invertible for every real \(z>\lambda_*\).  Hence, by the analytic
        Fredholm theorem (see, e.g.,~\cite[Theorem~VII.1.9]{Kato1995Book}), the non-invertible
        points of \(F(z)\) have no accumulation point in
        \(\{\operatorname{Re} z>\omega(a,V)\}\).

        On the other hand, \(F(\lambda_*)\) is not invertible, since
        \(P_{\lambda_*}\varphi=\varphi\) implies
        \((I-R_{\lambda_*}\mathcal{K})\varphi=0\).  Therefore there exists
        \(0<\delta<\delta_0\) such that \(F(z)\) is invertible whenever
        \(0<|z-\lambda_*|<\delta\).  By \eqref{eq:resolvent-factorization}, \(z-\mathcal{L}\) is
        also invertible for such \(z\).  Hence
        \[
        z\notin\sigma(\mathcal{L}),\qquad 0<|z-\lambda_*|<\delta,
        \]
        and so \(\lambda_*\) is isolated in \(\sigma(\mathcal{L})\).

        The geometric eigenspace of \(\mathcal{L}\) at \(\lambda_*\) is one-dimensional.  Indeed, if \(u\in D(\mathcal{A})\subset C(M;\mathbb{C})\) and \(\mathcal{L} u=\lambda_*u\), then Lemma~\ref{lem:eigen-resolvent} gives \(P_{\lambda_*}u=u\).  Since the eigenspace of \(P_{\lambda_*}\) for the eigenvalue one is spanned by \(\varphi\) (see, e.g.,~\cite[Lemma 3.1]{Liang2007CPAM}), the eigenspace of \(\mathcal{L}\) at \(\lambda_*\) is also spanned by \(\varphi\).

        It remains to rule out generalized eigenvectors.  We first claim that
        \(\varphi^*\) is strictly positive in the sense that
        \(\langle\varphi^*,h\rangle>0\) for every \(0\ne h\in C(M)_+\). Indeed,
        for such an \(h\), one has \(\int_M h\,\mathrm{d} m>0\).  Lemma~\ref{lem:minorization}
        gives \(P_{\lambda_*}^N h\ge c_N(\int_M h\,\mathrm{d} m)\mathbf{1}\).  Since
        \(\varphi^*\ge0\), \(\varphi^*\ne0\), we have \(\langle\varphi^*,\mathbf{1}\rangle>0\). Therefore, using
        \((P_{\lambda_*}^*)^N\varphi^*=\varphi^*\), we obtain
        \[
        \langle\varphi^*,h\rangle
        =\langle\varphi^*,P_{\lambda_*}^N h\rangle
        \ge c_N\Bigl(\int_M h\,\mathrm{d} m\Bigr)\langle\varphi^*,\mathbf{1}\rangle>0.
        \]

        Now suppose, to the contrary, that there exists
        \(u\in D(\mathcal{A})\subset C(M;\mathbb{C})\) such that
        \((\mathcal{L}-\lambda_*)u=\varphi\). Then
        \((\lambda_*-\mathcal{A})u-\mathcal{K} u=-\varphi\). Applying \(R_{\lambda_*}\) gives
        \((I-P_{\lambda_*})u=-R_{\lambda_*}\varphi\). Pairing this identity with
        \(\varphi^*\), and using \(P_{\lambda_*}^*\varphi^*=\varphi^*\), we obtain
        \[
        0
        =\langle\varphi^*,(I-P_{\lambda_*})u\rangle
        =-\langle\varphi^*,R_{\lambda_*}\varphi\rangle .
        \]
        This is impossible.  Indeed, \(\varphi\gg0\) and \eqref{eq:resolvent-formula} imply
        \(R_{\lambda_*}\varphi\gg0\), while the strict positivity of \(\varphi^*\)
        proved above gives \(\langle\varphi^*,R_{\lambda_*}\varphi\rangle>0\).
        Hence no generalized eigenvector exists.

        \medskip
        \textup{(iv)}
        By \textup{(ii)}, no spectral point of \(\mathcal{L}\) lies in
        \(\{\operatorname{Re} z>\lambda_*\}\).  It remains to exclude spectral points on the
        line \(\operatorname{Re} z=\lambda_*\).  Suppose that
        \(z=\lambda_*+i\beta\in\sigma(\mathcal{L})\).  By
        \eqref{eq:resolvent-factorization}, and since \(R_z\mathcal{K}\) is compact,
        the Fredholm
        alternative yields a nonzero \(u\in C(M;\mathbb{C})\) such that
        \(R_z\mathcal{K} u=u\).

        We first prove that \(u\) has a constant complex phase.  The modulus
        estimate gives
        \[
        |u|=|R_z\mathcal{K} u|
        \le R_{\lambda_*}|\mathcal{K} u|
        \le R_{\lambda_*}\mathcal{K}|u|
        =P_{\lambda_*}|u|.
        \]
        Set \(h:=P_{\lambda_*}|u|-|u|\). Then \(h\ge0\).  Moreover, using
        \(P_{\lambda_*}^*\varphi^*=\varphi^*\), we have
        \[
        \langle\varphi^*,h\rangle
        =
        \langle\varphi^*,P_{\lambda_*}|u|-|u|\rangle
        =
        \langle P_{\lambda_*}^*\varphi^*,|u|\rangle
        -
        \langle\varphi^*,|u|\rangle
        =0.
        \]
        Since \(\varphi^*\) is strictly positive, this is possible only if
        \(h=0\).  Hence \(P_{\lambda_*}|u|=|u|\).
        Thus, by \cite[Lemma 3.1]{Liang2007CPAM} again, \(|u|=c\varphi\) for some \(c>0\).
        Moreover, equality holds throughout the above chain of inequalities.  In
        particular,
        \(
        R_{\lambda_*}|\mathcal{K} u|=R_{\lambda_*}\mathcal{K}|u|.
        \)
        Since \(R_{\lambda_*}\) is injective, it follows that
        \(
        |\mathcal{K} u|=\mathcal{K}|u|.
        \)

        Set \(\zeta=u/|u|\).  Then \(\zeta\in C(M;S^1)\).  Fix
        \(x\in M\), and set \(\,\mathrm{d}\mu_x(y)=K(x,y)|u(y)|\,\mathrm{d} m(y)\).  Since
        \(|\mathcal{K} u(x)|=\mathcal{K}|u|(x)\), we have
        \[
        \left|\int_M\zeta(y)\,\mathrm{d}\mu_x(y)\right|
        =
        \int_M\,\mathrm{d}\mu_x(y).
        \]
        Choose \(\theta_x\in\mathbb{R}\) such that
        \(\int_M\zeta\,\mathrm{d}\mu_x=e^{i\theta_x}\mu_x(M)\).  Then
        \[
        \mu_x(M)
        =
        \int_M \operatorname{Re}\big(e^{-i\theta_x}\zeta(y)\big)\,\mathrm{d}\mu_x(y).
        \]
        Since \(\operatorname{Re}(e^{-i\theta_x}\zeta(y))\le1\), the last equality forces
        \(e^{-i\theta_x}\zeta(y)=1\) for \(\mu_x\)-a.e. \(y\).  For every
        \(x\in M\), assumption \textup{(K)} gives \(K(x,y)\ge\kappa_0>0\) for
        \(y\in B_g(x,\rho_0)\).  Since \(|u|=c\varphi\gg0\), the measure \(\mu_x\)
        has strictly positive density on \(B_g(x,\rho_0)\).  Hence
        \(\zeta=e^{i\theta_x}\) a.e. on \(B_g(x,\rho_0)\), and by continuity this
        identity holds everywhere on \(B_g(x,\rho_0)\).  Therefore \(\zeta\) is
        locally constant.  Since \(M\) is connected, \(\zeta\) is constant on
        \(M\).  Hence \(u=e^{i\theta_0}c\varphi\) for some
        \(\theta_0\in\mathbb{R}\).

        We next show that \(\beta=0\).  Since \(R_z\mathcal{K} u=u\) and
        \(u=e^{i\theta_0}c\varphi\), we have
        \(R_{\lambda_*+i\beta}\mathcal{K}\varphi=\varphi\).  On the other hand,
        \(P_{\lambda_*}\varphi=R_{\lambda_*}\mathcal{K}\varphi=\varphi\).  Thanks to \eqref{eq:resolvent-formula} again, for every \(x\in M\),
        \[
        \int_0^\infty e^{-i\beta t}\,\mathrm{d}\alpha_x(t)
        =
        \int_0^\infty \,\mathrm{d}\alpha_x(t),
        \]
        where
        \[
        \,\mathrm{d}\alpha_x(t)
        =
        e^{-\lambda_*t+\int_0^t a(\Phi_sx)\,\mathrm{d} s}
        (\mathcal{K}\varphi)(\Phi_tx)\,\mathrm{d} t .
        \]
        The density of \(\alpha_x\) is strictly positive and continuous.  Taking
        real parts gives
        \[
        \int_0^\infty \bigl(1-\cos(\beta t)\bigr)\,\mathrm{d}\alpha_x(t)=0.
        \]
        Since \(1-\cos(\beta t)\ge0\) and \(\alpha_x\) has strictly positive
        density, we obtain \(\cos(\beta t)=1\) for all \(t\ge0\).  Hence
        \(\beta=0\).  Therefore the line \(\operatorname{Re} z=\lambda_*\) contains no spectral
        point of \(\mathcal{L}\) except \(\lambda_*\), proving
        \eqref{eq:rightmost-spectrum}.
    \end{proof}

    Next, we present a sufficient condition to verify the condition given in the above theorem.
    \begin{proposition}\label{prop:divergence-condition}
        If
        \begin{equation}\label{eq:resolvent-mass-divergence}
            I(\lambda):=\int_M R_\lambda\mathbf{1}\,\mathrm{d} m\longrightarrow+\infty
            \qquad\text{as }\lambda\downarrow\omega(a,V),
        \end{equation}
        then
        \begin{equation}\label{eq:spectral-radius-divergence}
            {r}(P_\lambda)\longrightarrow+\infty
            \qquad\text{as }\lambda\downarrow\omega(a,V).
        \end{equation}
        In particular, \eqref{eq:spectral-radius-condition} holds for some
        \(\lambda_0>\omega(a,V)\), and the conclusions of
        Theorem~\ref{thm:spectral-radius} hold.
    \end{proposition}

    \begin{proof}
        Fix \(\Lambda>\omega(a,V)\).  By \eqref{eq:K-one-lower}, for
        \(\lambda\in(\omega(a,V),\Lambda]\),
        \[
        P_\lambda\mathbf{1}=R_\lambda\mathcal{K}\mathbf{1}\ge \kappa_M R_\lambda\mathbf{1}.
        \]
        Using Lemma~\ref{lem:minorization},
        \[
        P_\lambda^{N+1}\mathbf{1}(x)
        \ge c_N\kappa_M\int_M R_\lambda\mathbf{1}\,\mathrm{d} m
        =c_N\kappa_M I(\lambda).
        \]
        By positivity, for every \(n\ge1\),
        \[
        (P_\lambda^{N+1})^n\mathbf{1}\ge
        (c_N\kappa_M I(\lambda))^n\mathbf{1}.
        \]
        Hence \(\|(P_\lambda^{N+1})^n\|\ge
        (c_N\kappa_M I(\lambda))^n\).  By Gelfand's formula,
        \begin{equation}\label{eq:spectral-radius-lower-bound}
            {r}(P_\lambda)^{N+1}
            ={r}(P_\lambda^{N+1})
            \ge c_N\kappa_M I(\lambda).
        \end{equation}
        Thus \eqref{eq:resolvent-mass-divergence} and
        \eqref{eq:spectral-radius-lower-bound} imply
        \eqref{eq:spectral-radius-divergence}.  In particular,
        \eqref{eq:spectral-radius-condition} holds for \(\lambda\) close to
        \(\omega(a,V)\).
    \end{proof}

    \begin{remark}\label{rem:global-positive}
        If one assumes the stronger condition $K(x,y)\ge\kappa>0$ for all $x,y\in M$, then Lemma~\ref{lem:minorization} is unnecessary.  Indeed, $R_\lambda\mathcal{K} \mathbf{1}\ge \kappa I(\lambda)\mathbf{1}$, so the divergence of $I(\lambda)$ directly forces the spectral radius of $R_\lambda\mathcal{K} $ to diverge, which further implies that \eqref{eq:spectral-radius-condition} holds.
    \end{remark}

    \begin{theorem}\label{thm:CW-spectral-bound}
        For the operator \(\mathcal{L}=\mathcal{A}+\mathcal{K}\), the following statements are valid.
        \begin{enumerate}[label=\textup{(\roman*)},leftmargin=2.6em]
            \item Let \(\psi\in D(\mathcal{A})\), \(\psi\gg0\).  If \(\mathcal{L}\psi\ge\beta\psi\) for some \(\beta\in\mathbb{R}\), then \(\beta\le s(\mathcal{L})\).  If \(\mathcal{L}\psi\le\alpha\psi\) for some \(\alpha\in\mathbb{R}\), then \(s(\mathcal{L})\le\alpha\).

            \item Consequently, the spectral bound satisfies
            \begin{equation}\label{eq:CW-spectral-bound}
                \sup_{\psi\in D(\mathcal{A}),\,\psi\gg0}
                \inf_{x\in M}\frac{\mathcal{L}\psi(x)}{\psi(x)}
                \le
                s(\mathcal{L})
                \le
                \inf_{\psi\in D(\mathcal{A}),\,\psi\gg0}
                \sup_{x\in M}\frac{\mathcal{L}\psi(x)}{\psi(x)}.
            \end{equation}

            \item If, in addition, \eqref{eq:spectral-radius-condition} holds, then both inequalities in \eqref{eq:CW-spectral-bound} are equalities.  Their common value is the principal eigenvalue \(\lambda_*=s(\mathcal{L})\) obtained in Theorem~\ref{thm:spectral-radius}.
        \end{enumerate}
    \end{theorem}

    \begin{proof}
        The operator \(\mathcal{L}\) is resolvent positive because it is a bounded positive perturbation of the resolvent-positive transport generator \(\mathcal{A}\).  We shall use the following consequence of Thieme's resolvent-positive theory: for every real \(\lambda>s(\mathcal{L})\),
        \[
        (\lambda-\mathcal{L})^{-1}\ge0,
        \qquad
        {r}((\lambda-\mathcal{L})^{-1})
        =\frac1{\lambda-s(\mathcal{L})}.
        \]
        This follows from \cite[Proposition~3.4, Theorem~3.5, and Corollary~3.6]{Thieme1998DCDS}.  Consequently,
        ${r}((\lambda-\mathcal{L})^{-1})\to+
        \infty$ as \(\lambda\downarrow s(\mathcal{L})\).

        We first derive \textup{(i)} and \textup{(ii)} from this resolvent characterization.  Fix \(\psi\in D(\mathcal{A})\), \(\psi\gg0\).  Assume that \(\mathcal{L}\psi\le\alpha\psi\).  Suppose, to the contrary, that \(s(\mathcal{L})>\alpha\).  Then, for every \(\lambda>s(\mathcal{L})\),
        \[
        (\lambda-\mathcal{L})\psi
        \ge(\lambda-\alpha)\psi .
        \]
        Applying the positive resolvent gives
        \[
        0\le(\lambda-\mathcal{L})^{-1}\psi
        \le\frac1{\lambda-\alpha}\psi .
        \]
        Since \(\psi_{\min}:=\min_M\psi>0\), for every
        \(h\in C(M;\mathbb{C})\) we have
        \(|h|\le \psi_{\min}^{-1}\|h\|_\infty\psi\).  Using positivity once more,
        \[
        |(\lambda-\mathcal{L})^{-1}h|
        \le(\lambda-\mathcal{L})^{-1}|h|
        \le \frac{\|h\|_\infty}{\psi_{\min}(\lambda-\alpha)}\psi.
        \]
        Therefore
        \[
        \|(\lambda-\mathcal{L})^{-1}h\|_\infty
        \le \frac{\|\psi\|_\infty}{\psi_{\min}(\lambda-\alpha)}\|h\|_\infty .
        \]
        Thus the operator norms \(\|(\lambda-\mathcal{L})^{-1}\|\) remain bounded as \(\lambda\downarrow s(\mathcal{L})\), because \(s(\mathcal{L})>\alpha\).  This contradicts ${r}((\lambda-\mathcal{L})^{-1}) \rightarrow +\infty$ as \(\lambda\downarrow s(\mathcal{L})\).  Therefore \(s(\mathcal{L})\le\alpha\).

        The lower comparison is similar and shorter.  Assume that \(\mathcal{L}\psi\ge\beta\psi\).  If \(\beta>s(\mathcal{L})\), then \((\beta-\mathcal{L})^{-1}\ge0\), and
        \(
        (\beta-\mathcal{L})\psi\le0.
        \)
        Applying the positive resolvent gives \(\psi\le0\), contradicting \(\psi\gg0\).  Thus \(\beta\le s(\mathcal{L})\).  This proves \textup{(i)}.

        Now set
        \[
        \alpha_\psi:=\sup_{x\in M}\frac{\mathcal{L}\psi(x)}{\psi(x)},
        \qquad
        \beta_\psi:=\inf_{x\in M}\frac{\mathcal{L}\psi(x)}{\psi(x)}.
        \]
        Then \(\mathcal{L}\psi\le\alpha_\psi\psi\) and \(\mathcal{L}\psi\ge\beta_\psi\psi\).  By \textup{(i)},
        \[
        \beta_\psi\le s(\mathcal{L})\le\alpha_\psi .
        \]
        Taking the supremum over \(\psi\) on the left and the infimum over \(\psi\) on the right yields \eqref{eq:CW-spectral-bound}, proving \textup{(ii)}.

        It remains to prove \textup{(iii)}.  If condition \eqref{eq:spectral-radius-condition} holds, Theorem~\ref{thm:spectral-radius} gives a function \(\varphi\in D(\mathcal{A})\), \(\varphi\gg0\), such that
        \[
        \mathcal{L}\varphi=s(\mathcal{L})\varphi.
        \]
        For this test function the supremum and the infimum in \eqref{eq:CW-spectral-bound} are both equal to \(s(\mathcal{L})\).  Combined with \textup{(ii)}, this yields equality throughout.
    \end{proof}

    \begin{theorem}
        \label{thm:strict-averaged-threshold}
        For $\mathcal{L}=\mathcal{A}+\mathcal{K}$, the following statements are equivalent:
        \begin{enumerate}[label=\textup{(\roman*)},leftmargin=2.6em]
            \item $\mathcal{L}$ has a principal eigenvalue;
            \item $s(\mathcal{L})>\omega(a,V)$;
            \item there exists $\lambda_0>\omega(a,V)$ such that
            \(
            {r}(P_{\lambda_0})>1.
            \)
        \end{enumerate}
        Whenever these conditions hold, there is a unique
        $\lambda_*>\omega(a,V)$ satisfying \eqref{eq:critical-spectral-radius},
        and $\lambda_*=s(\mathcal{L})$.  Moreover, $\lambda_*$ is isolated and
        algebraically simple, has a strongly positive eigenfunction, and
        satisfies \eqref{eq:rightmost-spectrum}.
    \end{theorem}

    \begin{proof}
        \textup{(i)}$\Rightarrow$\textup{(ii)}.
        Let $\varphi\in D(\mathcal{A})$, $\varphi\gg0$, satisfy
        \[
        \mathcal{L}\varphi=s(\mathcal{L})\varphi.
        \]
        By \eqref{eq:omega} and the standard existence of
        maximizing invariant measures for continuous flows on compact spaces
        (see, e.g.,~\cite[Section~1]{MorroSantAnnaVarandas2020AHP}), choose
        $\nu_*\in\mathcal M_\Phi$ such that
        \[
        \int_M a\,\mathrm{d}\nu_*=\omega(a,V).
        \]
        Since $\varphi\gg0$, the chain rule for the generator gives
        $\log\varphi\in D(\mathcal{A})$ and
        \[
        V(\log\varphi)=\frac{V\varphi}{\varphi}.
        \]
        For every $t\ne0$, the $\Phi_t$-invariance of $\nu_*$ gives
        \[
        \int_M
        \frac{\log\varphi(\Phi_t x)-\log\varphi(x)}{t}
        \,\mathrm{d}\nu_*(x)=0.
        \]
        The integrand converges in $C(M)$ to $V(\log\varphi)$ as $t\to0$.
        Hence
        \[
        \int_M\frac{V\varphi}{\varphi}\,\mathrm{d}\nu_*
        =\int_M V(\log\varphi)\,\mathrm{d}\nu_*=0.
        \]
        Dividing the eigenvalue equation by $\varphi$ and integrating against
        $\nu_*$, we obtain the gap identity
        \begin{equation}\label{eq:principal-eigenvalue-gap}
            s(\mathcal{L})
            =\omega(a,V)+\int_M\frac{\mathcal{K}\varphi}{\varphi}\,\mathrm{d}\nu_*.
        \end{equation}
        Moreover, \eqref{eq:K-one-lower} gives
        \[
        \mathcal{K}\varphi
        \ge (\min_M\varphi)\mathcal{K}\mathbf{1}
        \ge \kappa_M\min_M\varphi>0
        \qquad\text{on }M.
        \]
        Thus the integral in \eqref{eq:principal-eigenvalue-gap} is strictly
        positive, proving \textup{(ii)}.

        \medskip
        \textup{(ii)}$\Rightarrow$\textup{(iii)}.
        Set
        \[
        \lambda^\dagger
        :=\frac{s(\mathcal{L})+\omega(a,V)}{2}>\omega(a,V).
        \]
        Suppose, to the contrary, that ${r}(P_{\lambda^\dagger})\le1$.
        If $z\in\mathbb{C}$ satisfies $\operatorname{Re} z>\lambda^\dagger$, then
        Lemma~\ref{lem:monotone-spectral-radius} gives
        \[
        {r}(P_{\operatorname{Re} z})<{r}(P_{\lambda^\dagger})\le1.
        \]
        The domination estimate \eqref{eq:complex-resolvent-domination} gives
        ${r}(R_z\mathcal{K})\le{r}(P_{\operatorname{Re} z})<1$.  Hence $I-R_z\mathcal{K}$ is
        invertible, and the factorization \eqref{eq:resolvent-factorization}
        shows that $z\in\rho(\mathcal{L})$.  Thus
        $\{z\in\mathbb{C}:\operatorname{Re} z>\lambda^\dagger\}\subset\rho(\mathcal{L})$, which implies
        $s(\mathcal{L})\le\lambda^\dagger$, a contradiction.  Consequently,
        ${r}(P_{\lambda^\dagger})>1$.

        \medskip
        \textup{(iii)}$\Rightarrow$\textup{(i)} follows directly from
        Theorem~\ref{thm:spectral-radius}.
    \end{proof}

    \begin{remark}\label{rem:strict-threshold-zero-transport}
        If $V=0$, then the growth bound in the absence of transport is
        $\omega(a,0)=\max_M a$, and
        Theorem~\ref{thm:strict-averaged-threshold} becomes the classical
        spectral-gap criterion
        \[
        \mathcal{L}\text{ has a principal eigenvalue}
        \quad\Longleftrightarrow\quad
        s(\mathcal{L})>\max_M a.
        \]
        This criterion is standard in the literature
        (see, e.g.,~\cite[Proposition~3.2]{ShenZhang2010JDE},
        \cite[Theorem~2.1]{BaoShen2017PAMS}, and
        \cite[Lemma~A.1(ii)]{WangZhang2026JMPA}).
    \end{remark}

    \begin{corollary}\label{cor:critical-coupling}
        For $\varepsilon>0$, set
        $\mathcal{L}_\varepsilon:=\mathcal{A}+\varepsilon\mathcal{K}$, and define
        \begin{equation}\label{eq:critical-coupling}
            r_*:=\sup_{\lambda>\omega(a,V)}r(P_\lambda)
            =\lim_{\lambda\downarrow\omega(a,V)}r(P_\lambda),
            \qquad
            \varepsilon_c:=\frac{1}{r_*},
        \end{equation}
        where $r_*\in(0,+\infty]$ and $1/(+\infty):=0$. Then $\mathcal{L}_\varepsilon$ admits a principal eigenvalue if and only if $\varepsilon>\varepsilon_c$. 
        For $0<\varepsilon\le\varepsilon_c$, no principal eigenvalue exists and $s(\mathcal{L}_\varepsilon)=\omega(a,V)$.
    \end{corollary}

    \begin{proof}
        The compact positive family associated with
        $\mathcal{L}_\varepsilon$ is
        \[
        R_\lambda(\varepsilon\mathcal{K})
        =\varepsilon P_\lambda.
        \]
        Hence Theorem~\ref{thm:strict-averaged-threshold} gives
        \[
        \mathcal{L}_\varepsilon\text{ has a principal eigenvalue}
        \quad\Longleftrightarrow\quad
        \varepsilon r(P_{\lambda_0})>1
        \text{ for some }\lambda_0>\omega(a,V).
        \]
        By the monotonicity of $\lambda\mapsto r(P_\lambda)$, this is
        equivalent to $\varepsilon r_*>1$, or
        $\varepsilon>\varepsilon_c$.

        Since
        $\mathcal{L}_\varepsilon$ has no principal eigenvalue,
        Theorem~\ref{thm:strict-averaged-threshold} gives
        $s(\mathcal{L}_\varepsilon)\le\omega(a,V)$.  On the other hand,
        monotonicity of the spectral bound under positive bounded perturbations
        yields
        \[
        s(\mathcal{L}_\varepsilon)
        \ge s(\mathcal{A})=\omega(a,V);
        \]
        see \cite[Theorem 4.2]{Thieme1998DCDS}.  The conclusion follows.
    \end{proof}

    \begin{remark}\label{rem:critical-coupling-bounds}
        Notice that
        \(
        r(P_\lambda)=r(\mathcal{K}R_\lambda)
        \le \|\mathcal{K}R_\lambda\|.
        \)
        Consequently, if
        \[
        M_{\mathcal{K}}
        :=\sup_{\lambda>\omega(a,V)}
        \|\mathcal{K}R_\lambda\|<+\infty,
        \]
        then $\varepsilon_c\ge M_{\mathcal{K}}^{-1}>0$, and
        $\mathcal{L}_\varepsilon$ has no principal eigenvalue whenever
        $0<\varepsilon\le M_{\mathcal{K}}^{-1}$.

        A convenient sufficient condition is
        \[
        \sup_{\lambda>\omega(a,V)}
        \int_M R_\lambda\mathbf{1}\,\mathrm{d}m<+\infty.
        \]
        Indeed, positivity gives
        \[
        \|\mathcal{K}R_\lambda\|
        \le \|K\|_\infty
        \int_M R_\lambda\mathbf{1}\,\mathrm{d}m.
        \]
        In the opposite direction, if the resolvent mass in
        \eqref{eq:resolvent-mass-divergence} diverges, then
        Proposition~\ref{prop:divergence-condition} gives $r_*=+\infty$ and
        hence $\varepsilon_c=0$.
    \end{remark}

    \section{Uniquely ergodic flows}\label{sec:uniquely-ergodic}

    We next apply the criterion of Section~\ref{sec:abstract-criterion} to uniquely ergodic flows. Throughout this section, we keep the standing assumptions and notation of Section~\ref{sec:abstract-criterion} and assume in addition that $\Phi_t$ is uniquely ergodic.  Let \(\mu\) denote its unique invariant Borel probability measure, and set
    \[
    \bar a_\mu:=\int_M a\,\mathrm{d}\mu.
    \]

    The invariant measure of a smooth minimal uniquely ergodic flow may be singular
     with respect to Riemannian volume, rather than absolutely continuous. Such examples can
    be obtained by suspending smooth circle diffeomorphisms with singular unique
    invariant measure (see, e.g.,~\cite{Sadovskaya2009ETDS}). Thus the
    volume-compatibility used in Theorem~\ref{thm:intro-B} is not automatic.

    \begin{lemma}\label{lem:uniform-ergodic}
        For every \(b\in C(M)\),
        \begin{equation}\label{eq:uniform-ergodic-average}
            \frac{1}{t}\int_0^t b(\Phi_sx)\,\mathrm{d}s
            \longrightarrow \int_M b\,\mathrm{d}\mu
            \qquad\text{uniformly for }x\in M.
        \end{equation}
        Consequently,
        \begin{equation}\label{eq:UE-threshold}
            \omega(a,V)=\bar a_\mu=\int_M a\,\mathrm{d}\mu.
        \end{equation}
    \end{lemma}

    \begin{proof}
        The convergence in \eqref{eq:uniform-ergodic-average} is the
        continuous-time version of the uniform ergodic theorem for
        uniquely ergodic systems; compare \cite[Theorem~6.19]{Walters1982Book} for
        discrete time.  The standard proof, using weak-* limits of the occupation
        measures, applies verbatim to flows.  Applying \eqref{eq:uniform-ergodic-average}
        with \(b=a\), then taking the supremum in \(x\), gives
        \eqref{eq:UE-threshold}.
    \end{proof}

    We now prove Theorem~\ref{thm:intro-B}.
    \begin{proof}[Proof of Theorem~\ref{thm:intro-B}]
        By \eqref{eq:UE-threshold},
        $\omega(a,V)=\bar a_\mu$.  Hence, by
        Theorem~\ref{thm:spectral-radius} and
        Proposition~\ref{prop:divergence-condition}, it is enough to verify
        the divergence condition \eqref{eq:resolvent-mass-divergence}, namely
        \[
        I(\lambda):=\int_M R_\lambda\mathbf{1}\,\mathrm{d} m\longrightarrow +\infty
        \qquad\text{as }\lambda\downarrow\bar a_\mu.
        \]

        First assume \textup{(BD)}.  From the resolvent formula \eqref{eq:resolvent-formula} and assumption \textup{(BD)}, for every $\lambda>\bar a_\mu$ and every $x\in M$,
        \begin{align*}
            R_\lambda\mathbf{1}(x)
            &=\int_0^\infty
            \exp\left(-\lambda t+
            \int_0^t a(\Phi_sx)\,\mathrm{d} s\right)\,\mathrm{d} t \\
            &\ge \int_0^\infty
            \exp\left(-\lambda t+\bar a_\mu t-C_0\right)\,\mathrm{d} t
            =\frac{e^{-C_0}}{\lambda-\bar a_\mu}.
        \end{align*}
        Therefore
        \[
        I(\lambda)\ge \frac{e^{-C_0}}{\lambda-\bar a_\mu}\to+\infty
        \quad \text{ as }\lambda\downarrow\bar a_\mu.
        \]

        Next assume \textup{(AC)}.  Put
        \[
        B_t(x)=\int_0^t\bigl(a(\Phi_sx)-\bar a_\mu\bigr)\,\mathrm{d} s,
        \qquad t\ge0,\quad x\in M.
        \]
        Since $\mu$ is invariant under $\Phi_t$, Fubini's theorem (see, e.g.,~\cite[Theorem~2.37]{Folland1999Book}) gives
        \[
        \int_M B_t(x)\,\mathrm{d}\mu(x)
        =\int_0^t\left(\int_M a(\Phi_sx)\,\mathrm{d}\mu(x)-\bar a_\mu\right)\,\mathrm{d} s=0.
        \]
        Jensen's inequality for the probability measure \(\mu\) (see, e.g.,~\cite[Theorem~3.3]{Rudin1987Book}) gives, for every $t\ge0$,
        \[
        \int_M e^{B_t(x)}\,\mathrm{d}\mu(x)
        \ge \exp\left(\int_M B_t\,\mathrm{d}\mu\right)=1.
        \]
        Because $\,\mathrm{d}\mu=\rho\,\mathrm{d} m$ and $\rho\le C_1$, every nonnegative measurable function $h$ satisfies
        \[
        \int_M h\,\mathrm{d} m\ge \frac1{C_1}\int_M h\,\mathrm{d}\mu.
        \]
        Using Tonelli's theorem (see, e.g.,~\cite[Theorem~2.37]{Folland1999Book}) and the preceding Jensen bound, we obtain
        \begin{align*}
            I(\lambda)
            &=\int_0^\infty e^{-(\lambda-\bar a_\mu)t}
            \int_M e^{B_t(x)}\,\mathrm{d} m(x)\,\mathrm{d} t \\
            &\ge \frac1{C_1}\int_0^\infty e^{-(\lambda-\bar a_\mu)t}
            \int_M e^{B_t(x)}\,\mathrm{d}\mu(x)\,\mathrm{d} t \\
            &\ge \frac1{C_1}\int_0^\infty e^{-(\lambda-\bar a_\mu)t}\,\mathrm{d} t
            =\frac1{C_1(\lambda-\bar a_\mu)}.
        \end{align*}
        Thus the required divergence holds also under \textup{(AC)}.
    \end{proof}

    We spell out the effect of a positive time
    change.
    Let $f\in C^1(M)$ with $f>0$, and set
    \[
    W=fV.
    \]
    The standard time-change formula yields the following proposition.
    \begin{proposition}\label{prop:time-change-property}
        If the $V$-flow is uniquely ergodic with invariant probability measure $\mu$, then the $W$-flow is uniquely ergodic, and its unique invariant probability measure is
        \begin{equation}\label{eq:mu-f}
            \,\mathrm{d}\mu_f=\frac{f^{-1}}{\int_M f^{-1}\,\mathrm{d}\mu}\,\mathrm{d}\mu.
        \end{equation}
        Consequently, for every $a\in C(M)$,
        \begin{equation}\label{eq:time-change-average}
            \omega(a,W)=\int_M a\,\mathrm{d}\mu_f
            =\frac{\int_M a f^{-1}\,\mathrm{d}\mu}{\int_M f^{-1}\,\mathrm{d}\mu}.
        \end{equation}
    \end{proposition}

    Write \(\bar a_{\mu_f}:=\int_M a\,\mathrm{d}\mu_f\).  For the time-changed statement, we use the corresponding assumptions for the $W$-flow $\Psi_t$ and its invariant measure $\mu_f$:
    \begin{description}[leftmargin=3.2em]
        \item[(BD$_f$)] there exists $C_0>0$ such that, for all $x\in M$ and $t\ge0$,
        \[
        \int_0^t\bigl(a(\Psi_sx)-\bar a_{\mu_f}\bigr)\,\mathrm{d} s\ge -C_0;
        \]
        \item[(AC$_f$)] there exists $C_1>0$ such that
        \[
        \,\mathrm{d}\mu_f=\rho\,\mathrm{d} m,
        \qquad
        0\le\rho\le C_1\quad m\text{-a.e.}
        \]
    \end{description}
    When $f\equiv1$, these are exactly \textup{(BD)} and \textup{(AC)} from Theorem~\ref{thm:intro-B}.

    \begin{corollary}\label{cor:time-changed-UE}
        Suppose that the $V$-flow is uniquely ergodic with invariant probability measure $\mu$.  Let $W=fV$, where $f\in C^1(M)$ and $f>0$, let $\mu_f$ be given by \eqref{eq:mu-f}, and write $\bar a_{\mu_f}:=\int_M a\,\mathrm{d}\mu_f$.  If either \textup{(BD$_f$)} or \textup{(AC$_f$)} holds, then $W+a+\mathcal{K}$ has a principal eigenvalue
        \[
        \lambda_*=s(W+a+\mathcal{K})>\bar a_{\mu_f}.
        \]
    \end{corollary}

    \begin{proof}
        By the positive time-change property, the $W$-flow is uniquely ergodic with invariant probability measure $\mu_f$, and \eqref{eq:time-change-average} gives $\omega(a,W)=\bar a_{\mu_f}$.  Applying Theorem~\ref{thm:intro-B} with $W$ in place of $V$ gives the conclusion.
    \end{proof}

    We remark that Proposition~\ref{prop:time-change-property} and
    Corollary~\ref{cor:time-changed-UE} are included only to make explicit how the
    invariant measure and the averaged growth rate transform under the change
    \(V\mapsto W=fV\).

    \begin{remark}
        The results in this section apply to flows that are uniquely ergodic on the
        whole manifold.  Many natural transport fields, however, fail to have this
        global property, although their orbit closures still possess a simple compact
        geometric structure.  A basic example is a rotation on a sphere.  For instance,
        on \(S^2\),
        \[
        x'=-y,\qquad y'=x,\qquad z'=0,
        \]
        the two poles are fixed points, while every other orbit is a periodic circle
        of latitude.  Hence the flow is not uniquely ergodic on the
        whole sphere.  Such flows are not covered directly by
        Theorem~\ref{thm:intro-B}.

        In the next section, we treat such cases by decomposing the dynamics along
        compact orbit closures and reducing the remaining problem to the corresponding
        quotient space.
    \end{remark}
    \section{Compact abelian group actions and quotient reduction}\label{sec:group}

    In this section, we apply the criterion of Section~\ref{sec:abstract-criterion} to flows whose orbit closures are compact group orbits.  Throughout this section, we still use the standing assumptions and notation of Section~\ref{sec:abstract-criterion} and impose the following additional dynamical structure.

    Let \(H\) be a compact abelian Lie group acting smoothly and isometrically on \(M\).  The Riemannian volume \(m\) is then \(H\)-invariant.  Let \(\gamma:\mathbb{R}\to H\) be a smooth one-parameter subgroup with dense image, and assume that the transport field \(V\) is generated by \(\gamma\):
    \[
    V(x):=\frac{\mathrm{d}}{\mathrm{d}t}\bigg|_{t=0}\gamma(t)\cdot x,
    \qquad
    \Phi_t x:=\gamma(t)\cdot x.
    \]
    For \(x\in M\), set \(F_x:=H\cdot x\).  Since \(\gamma(\mathbb{R})\) is dense in \(H\), the sets \(F_x\) are exactly the orbit closures of the \(V\)-flow.  Let
    \[
    Q:=M/H,
    \qquad
    \pi:M\to Q
    \]
    be the orbit space and quotient map.  The space $Q$ is compact and metrizable, even though it may lack a smooth manifold structure and have leaves of varying dimensions

    Let \(\mathrm{d}h\) denote normalized Haar measure on \(H\).  For \(x\in M\), define the probability measure \(\nu_x\) on \(F_x\) by
    \begin{equation}\label{eq:Haar-leaf-measure}
        \int_{F_x}\psi\,\mathrm{d}\nu_x
        =\int_H\psi(h\cdot x)\,\mathrm{d}h,
        \qquad \psi\in C(M).
    \end{equation}
    If \(q=\pi(x)\), we write \(F_q=F_x\) and \(\nu_q=\nu_x\).  This is well-defined because replacing \(x\) by another point in the same \(H\)-orbit gives the same pushed-forward Haar measure.  The standard identification \(F_x\simeq H/H_x\) and the invariant measure on this compact homogeneous space are discussed in the literature (see, e.g.,~\cite[Chapter~I, Sections~3--4]{Bredon1972Book} and \cite[Proposition~2.46 and Theorem~2.51]{Folland2016Book}).

    Define the Haar leaf average of \(a\) by
    \begin{equation}\label{eq:Haar-leaf-average}
        \bar a(q)=\int_{F_q}a\,\mathrm{d}\nu_q
        =\int_H a(h\cdot x)\,\mathrm{d}h,
        \qquad q=\pi(x).
    \end{equation}
    The function \(\bar a\) is continuous on \(Q\), since the right-hand side of \eqref{eq:Haar-leaf-average} is a continuous \(H\)-invariant function of \(x\).

    \paragraph{Relation with ambient compact Lie groups.}
    This direct $H$-formulation includes the common situation where the flow is
    generated by a one-parameter subgroup of an ambient compact Lie group $G$.
    If $G$ acts smoothly and isometrically on $M$, $\xi$ belongs to the Lie
    algebra of $G$, and
    \[
    V_\xi(x)=\frac{\,\mathrm{d}}{\,\mathrm{d} t}\bigg|_{t=0}\exp_G(t\xi)\cdot x,
    \]
    then
    \[
    H=\overline{\{\exp_G(t\xi):t\in\mathbb{R}\}}
    \]
    is a compact abelian Lie subgroup of $G$, and the above framework holds
    with $\gamma(t)=\exp_G(t\xi)$.

    \begin{lemma}
        \label{lem:leafwise-Haar-threshold}
        For every $x\in M$, the restricted flow $\Phi_t|_{F_x}$ is uniquely
        ergodic, and its unique invariant probability measure is the
        $H$-invariant probability measure $\nu_x$ defined by
        \eqref{eq:Haar-leaf-measure}.
        Consequently, for the group-generated vector field $V$,
        \begin{equation}\label{eq:Haar-threshold}
            \omega(a,V)=\max_{q\in Q}\bar a(q).
        \end{equation}
    \end{lemma}

    \begin{proof}
        The general theory of Haar measures on homogeneous spaces gives a unique \(H\)-invariant probability measure on \(F_x\) (see, e.g.,~\cite[Proposition~2.46 and Theorem~2.51]{Folland2016Book}).  The dense-subgroup argument below shows that this measure is also the unique probability measure invariant under the restricted flow.

        Fix $x\in M$.  Recall that $\{\gamma(t):t\in\mathbb{R}\}$ is dense in $H$.
        The measure $\nu_x$ is $H$-invariant, hence invariant under
        $\Phi_t|_{F_x}$.

        Conversely, let $\mu$ be an arbitrary $\Phi_t|_{F_x}$-invariant
        probability measure.  For $h\in H$, choose a sequence
        $t_n\in\mathbb{R}$ such that $\gamma(t_n)\to h$ in $H$.
        For any $\psi\in C(F_x)$, the map
        $y\mapsto\psi(\gamma(t_n)\cdot y)$ converges uniformly to
        $y\mapsto\psi(h\cdot y)$; therefore, by the dominated convergence
        theorem and the $\Phi_t|_{F_x}$-invariance of $\mu$,
        \[
        \int_{F_x}\psi(h\cdot y)\,\mathrm{d}\mu(y)
        =\lim_{n\to\infty}
        \int_{F_x}\psi(\gamma(t_n)\cdot y)\,\mathrm{d}\mu(y)
        =\int_{F_x}\psi(y)\,\mathrm{d}\mu(y).
        \]
        Thus $\mu$ is $H$-invariant.  Averaging this identity over $H$ and
        applying Fubini's theorem gives
        \[
        \int_{F_x}\psi\,\mathrm{d}\mu
        =\int_{F_x}\Bigl(\int_H\psi(h\cdot y)\,\mathrm{d} h\Bigr)\,\mathrm{d}\mu(y).
        \]
        The inner Haar average is constant along the orbit $H\cdot x = F_x$ and
        equals $\int_H\psi(h\cdot x)\,\mathrm{d} h$.  Hence
        \[
        \int_{F_x}\psi\,\mathrm{d}\mu
        =\int_H\psi(h\cdot x)\,\mathrm{d} h
        =\int_{F_x}\psi\,\mathrm{d}\nu_x.
        \]
        Since $\psi\in C(F_x)$ was arbitrary, $\mu = \nu_x$.  As $\mu$ was an
        arbitrary invariant probability measure, the restricted flow
        $\Phi_t|_{F_x}$ is uniquely ergodic with unique invariant measure $\nu_x$.

        Finally, let \(\eta\) be any \(\Phi_t\)-invariant probability measure
        on \(M\).  By the same dense-subgroup argument used above, now applied
        on \(M\), the measure \(\eta\) is \(H\)-invariant.  Hence, by
        \(H\)-invariance and Fubini's theorem,
        \[
        \int_M a\,\mathrm{d}\eta
        =
        \int_M\int_H a(h\cdot x)\,\mathrm{d}h\,\mathrm{d}\eta(x)
        =
        \int_M \bar a(\pi(x))\,\mathrm{d}\eta(x)
        \leq \max_{q\in Q}\bar a(q).
        \]
        Conversely, each \(\nu_q\), regarded as a probability measure on \(M\)
        supported on \(F_q\), is \(\Phi_t\)-invariant and satisfies
        \[
        \bar{a}(q)=\int_M a \,\mathrm{d}\nu_q \le \max_{\nu\in\mathcal{M}_\Phi}\int_M a \,\mathrm{d}\nu.
        \]
        Therefore, by the invariant-measure characterization of
        \(\omega(a,V)\),
        \[
        \omega(a,V)=\max_{q\in Q}\bar a(q).
        \]
        The compactness of \(Q\) and the continuity of \(\bar a\) ensure that
        the maximum is attained.
    \end{proof}

    The following disintegration formula follows from Haar averaging,
    the $H$-invariance of $m$, and the definition $\eta=\pi_\#m$ (see, e.g.,~\cite[Chapter~II, Section~4]{Helgason2001Book}).
    \begin{lemma}\label{lem:Haar-disintegration}
        Let $\eta=\pi_\#m$.  Then, for every $\psi\in C(M)$,
        \begin{equation}\label{eq:Haar-disintegration}
            \int_M \psi(x)\,\mathrm{d} m(x)
            =\int_Q\int_{F_q}\psi\,\mathrm{d}\nu_q\,\mathrm{d}\eta(q).
        \end{equation}
        In particular, the above formula remains valid in the presence of singular orbits.
    \end{lemma}
    \begin{proof}
        For readers' convenience, we present a short proof.
        Define
        \[
        \Psi(x):=\int_H\psi(h\cdot x)\,\mathrm{d} h .
        \]
        Then \(\Psi\) is continuous and \(H\)-invariant, hence it descends to a
        continuous function on \(Q\), still denoted by \(\Psi(q)\).  By the
        definition of \(\nu_q\),
        \[
        \Psi(q)=\int_{F_q}\psi\,\mathrm{d}\nu_q .
        \]
        Since \(m\) is \(H\)-invariant, Fubini's theorem gives
        \[
        \int_M\Psi(x)\,\mathrm{d} m(x)
        =
        \int_H\int_M\psi(h\cdot x)\,\mathrm{d} m(x)\,\mathrm{d} h
        =
        \int_M\psi(x)\,\mathrm{d} m(x).
        \]
        On the other hand, since \(\eta=\pi_\#m\),
        \[
        \int_M\Psi(x)\,\mathrm{d} m(x)
        =
        \int_Q\Psi(q)\,\mathrm{d}\eta(q).
        \]
        Combining the last two identities yields
        \eqref{eq:Haar-disintegration}.
    \end{proof}

    We now prove Theorem~\ref{thm:intro-C}.

    \begin{proof}[Proof of Theorem~\ref{thm:intro-C}]
        Set
        \[
        \bar a_*:=\max_{q\in Q}\bar a(q).
        \]
        Identity \eqref{eq:Haar-threshold} gives
        $\omega(a,V)=\bar a_*$.  By
        Proposition~\ref{prop:divergence-condition} and
        Theorem~\ref{thm:spectral-radius}, it is enough to verify
        \eqref{eq:resolvent-mass-divergence}, namely
        \[
        I(\lambda):=\int_M R_\lambda\mathbf{1}\,\mathrm{d} m\to+\infty
        \qquad(\lambda\downarrow\bar a_*),
        \]
        where
        \[
        R_\lambda=(\lambda-V-a)^{-1}.
        \]

        Fix $q\in Q$.  For $x\in F_q$, set
        \[
        B_t(x)=\int_0^t\bigl(a(\Phi_sx)-\bar a(q)\bigr)\,\mathrm{d} s.
        \]
        Since $\nu_q$ is invariant under the $V$-flow on $F_q$ and
        $\int_{F_q}a\,\mathrm{d}\nu_q=\bar a(q)$, we have
        \begin{align*}
            \int_{F_q}B_t\,\mathrm{d}\nu_q
            &=\int_0^t\left(\int_{F_q}a(\Phi_sx)\,\mathrm{d}\nu_q(x)-\bar a(q)\right)\,\mathrm{d} s=0.
        \end{align*}
        Jensen's inequality (see, e.g.,~\cite[Theorem~3.3]{Rudin1987Book}) therefore gives the
        estimate
        \[
        \int_{F_q}e^{B_t}\,\mathrm{d}\nu_q
        \ge
        \exp\left(\int_{F_q}B_t\,\mathrm{d}\nu_q\right)
        =1,
        \qquad q\in Q,
        \quad t\ge0.
        \]

        Using the resolvent formula leaf by leaf, for $\lambda>\bar a_*$ we get
        \begin{equation}\label{eq:leafwise-resolvent-lower-bound}
            \begin{aligned}
                \int_{F_q}R_\lambda\mathbf{1}\,\mathrm{d}\nu_q
                &=\int_0^\infty e^{-(\lambda-\bar a(q))t}
                \int_{F_q}e^{B_t(x)}\,\mathrm{d}\nu_q(x)\,\mathrm{d} t \\
                &\ge \frac{1}{\lambda-\bar a(q)}.
            \end{aligned}
        \end{equation}
        Integrating \eqref{eq:leafwise-resolvent-lower-bound} by means of
        \eqref{eq:Haar-disintegration}, we obtain
        \begin{equation}\label{eq:quotient-resolvent-lower-bound}
            I(\lambda)
            \ge
            \int_Q\frac{\,\mathrm{d}\eta(q)}{\lambda-\bar a(q)}.
        \end{equation}
        The right-hand side of \eqref{eq:quotient-resolvent-lower-bound}
        diverges by \eqref{eq:quotient-integral-condition}.  Proposition~\ref{prop:divergence-condition}
        and Theorem~\ref{thm:spectral-radius} then imply the existence of a
        principal eigenvalue with spectral bound strictly above $\bar a_*$.
    \end{proof}

    \begin{remark}
        Write $\bar a_*:=\max_{q\in Q}\bar a(q)$.  The
        condition \eqref{eq:quotient-integral-condition} is automatic if the
        maximum set of $\bar a$ has positive
        $\eta$-measure.  More generally, the condition depends jointly on the
        behavior of $\bar a$ near its maximizers and on the local mass of the
        quotient measure $\eta$.

        For example, suppose that a neighborhood of a maximizer $q_0$ admits
        coordinates in $\mathbb{R}^k$ in which $\eta$ is comparable to
        $k$-dimensional Lebesgue measure, and that
        \[
        \bar a_* - \bar a(q)\le C|q-q_0|^\alpha
        \]
        near $q_0$.  Then the quotient integral diverges whenever
        $\alpha\ge k$.  In particular, a nondegenerate quadratic maximum gives
        divergence in quotient dimensions one and two.
    \end{remark}

    If \(W=fV\) with \(f>0\), then the \(W\)-flow has the same orbit closures as the \(V\)-flow.  On each orbit closure, the invariant probability measure $\nu_q$ is replaced by its inverse-speed reweighting, with density proportional to \(f^{-1}\).

    \begin{corollary}\label{cor:time-changed-group}
        Let \(f\in C^1(M)\) with \(f>0\), and set
        \[
        W=fV .
        \]
        For \(q\in Q\), define
        \begin{equation}\label{eq:time-changed-leaf-average}
            \bar a_f(q)
            =
            \frac{\int_{F_q}a f^{-1}\,\mathrm{d}\nu_q}
            {\int_{F_q}f^{-1}\,\mathrm{d}\nu_q},
            \qquad
            \bar a_{f,*}:=\max_{q\in Q}\bar a_f(q).
        \end{equation}
        Then \(\bar a_f\in C(Q)\) and
        \begin{equation}\label{eq:time-changed-group-threshold}
            \omega(a,W)=\bar a_{f,*}.
        \end{equation}
        Moreover, if
        \begin{equation}\label{eq:quotient-integral-condition-f}
            \int_Q \frac{\,\mathrm{d}\eta(q)}{\lambda-\bar a_f(q)}
            \longrightarrow +\infty
            \qquad\text{as }\lambda\downarrow\bar a_{f,*},
        \end{equation}
        then $W+a+\mathcal{K}$ has a principal eigenvalue
        \[
        \lambda_*=s(W+a+\mathcal{K})>\bar a_{f,*}.
        \]
    \end{corollary}

    \begin{proof}
        For each \(q\in Q\), the \(W\)-invariant probability measure on \(F_q\) is
        \begin{equation}\label{eq:time-changed-leaf-measure}
            \,\mathrm{d}\nu_q^f=
            \frac{f^{-1}\,\mathrm{d}\nu_q}{\int_{F_q}f^{-1}\,\mathrm{d}\nu_q}.
        \end{equation}
        By \eqref{eq:time-changed-leaf-average},
        \(\bar a_f(q)=\int_{F_q}a\,\mathrm{d}\nu_q^f\).
        Let $\Psi_t$ be the flow of $W=fV$.  For each $x\in M$ there is an
        increasing $C^1$ function $\tau=\tau(t,x)$ such that
        \[
        \Psi_t x=\Phi_{\tau(t,x)}x,
        \qquad
        \frac{\,\mathrm{d}\tau}{\,\mathrm{d} t}=f(\Phi_{\tau(t,x)}x).
        \]
        Equivalently,
        \[
        t=\int_0^{\tau(t,x)}f(\Phi_sx)^{-1}\,\mathrm{d} s.
        \]
        Thus $\Phi_t$ and $\Psi_t$ have the same orbits and the same orbit
        closures $F_q=H\cdot x$.

        Using this identity, for $t>0$,
        \begin{equation}\label{eq:time-change-ratio}
            \frac1t\int_0^t a(\Psi_sx)\,\mathrm{d} s
            =
            \frac{\int_0^{\tau(t,x)}a(\Phi_sx)f(\Phi_sx)^{-1}\,\mathrm{d} s}
            {\int_0^{\tau(t,x)}f(\Phi_sx)^{-1}\,\mathrm{d} s}.
        \end{equation}
        Since $f$ is positive and continuous on compact $M$, $\tau(t,x)\to\infty$
        as $t\to\infty$ uniformly in $x$.  Applying
        Lemma~\ref{lem:leafwise-Haar-threshold} to $a f^{-1}$ and to
        $f^{-1}$ gives uniform convergence of the ratio in
        \eqref{eq:time-change-ratio} to $\bar a_f(\pi(x))$.  This proves
        \eqref{eq:time-changed-group-threshold}.  The same time-change formula
        also shows that the measure in \eqref{eq:time-changed-leaf-measure} is
        the invariant probability measure of
        the $W$-flow on $F_q$.

        It remains to verify the divergence condition needed in the abstract
        spectral-radius criterion.  Let
        \[
        R_\lambda^W=(\lambda-W-a)^{-1}.
        \]
        For $x\in F_q$, set
        \[
        B_t^f(x)=\int_0^t\bigl(a(\Psi_sx)-\bar a_f(q)\bigr)\,\mathrm{d} s.
        \]
        Since $\nu_q^f$ is invariant under $\Psi_t$ and
        $\int_{F_q}a\,\mathrm{d}\nu_q^f=\bar a_f(q)$, we have
        \[
        \int_{F_q}B_t^f\,\mathrm{d}\nu_q^f=0.
        \]
        By Jensen's inequality for probability measures,
        \[
        \int_{F_q}e^{B_t^f}\,\mathrm{d}\nu_q^f
        \ge
        \exp\left(\int_{F_q}B_t^f\,\mathrm{d}\nu_q^f\right)=1.
        \]
        Since $f$ is positive and continuous on compact $M$, the measures
        $\nu_q^f$ and $\nu_q$ are uniformly equivalent.  Hence there exists
        a constant $c_f>0$, independent of $q$ and $t$, such that
        \[
        \int_{F_q}e^{B_t^f}\,\mathrm{d}\nu_q\ge c_f.
        \]
        Consequently, for $\lambda>\bar a_{f,*}$,
        \[
        \int_{F_q}R_\lambda^W\mathbf{1}\,\mathrm{d}\nu_q
        \ge
        c_f\int_0^\infty e^{-(\lambda-\bar a_f(q))t}\,\mathrm{d} t
        =
        \frac{c_f}{\lambda-\bar a_f(q)}.
        \]
        By \eqref{eq:Haar-disintegration},
        \begin{equation}\label{eq:time-changed-quotient-resolvent-lower-bound}
            \int_M R_\lambda^W\mathbf{1}\,\mathrm{d} m
            \ge
            c_f
            \int_Q\frac{\,\mathrm{d}\eta(q)}{\lambda-\bar a_f(q)}.
        \end{equation}
        The right-hand side of
        \eqref{eq:time-changed-quotient-resolvent-lower-bound} diverges by
        \eqref{eq:quotient-integral-condition-f}.  Proposition~\ref{prop:divergence-condition}
        and Theorem~\ref{thm:spectral-radius} give the desired principal
        eigenvalue.
    \end{proof}

    \section{Applications}\label{sec:applications}

    We now spell out the quotient criterion in concrete geometries.  The torus examples show how the dimension of the orbit closure determines the averaging scale.  The two-dimensional examples show that transport in a rational direction may either create or suppress a principal eigenvalue.  The sphere examples show that the same quotient integral remains meaningful when singular leaves occur.

    Throughout this section, for nonnegative functions $F$ and $G$, the notation
    $F\asymp G$ near a specified point means that there exist constants
    $0<c\le C<\infty$ such that $cG\le F\le CG$ in a neighborhood of that
    point.

    \subsection{Torus quotients}

    Let $M=\mathbb{T}^n=\mathbb{R}^n/\mathbb{Z}^n$ and let
    \[
    V_b=b\cdot\nabla,
    \qquad b\in\mathbb{R}^n\setminus\{0\}.
    \]
    Define the subtorus
    \[
    H_b=\overline{\{tb\bmod\mathbb{Z}^n:t\in\mathbb{R}\}}\subset\mathbb{T}^n.
    \]
    The leaves are
    \[
    F_x=x+H_b,
    \]
    and the quotient is
    \[
    Q_b=\mathbb{T}^n/H_b.
    \]
    Let $\pi_b:\mathbb{T}^n\to Q_b$ be the quotient map, let $m_{H_b}$ be the
    normalized Haar measure on $H_b$, and let
    $\eta_b:=(\pi_b)_\# m$ be the normalized Haar measure on $Q_b$.  For
    $q=x+H_b\in Q_b$, the leaf average and the corresponding transport growth
    bound are
    \begin{equation}\label{eq:torus-leaf-threshold}
        \bar a(q)=\int_{H_b}a(x+h)\,\mathrm{d} m_{H_b}(h),
        \qquad
        \omega(a,V_b)=\max_{q\in Q_b}\bar a(q).
    \end{equation}
    By Theorem~\ref{thm:intro-C}, with the transport growth bound in
    \eqref{eq:torus-leaf-threshold}, the operator $V_b+a+\mathcal{K}$ has a
    principal eigenvalue whenever
    \begin{equation}\label{eq:torus-integral-condition}
        \int_{Q_b}\frac{\,\mathrm{d}\eta_b(q)}{\lambda-\bar a(q)}\to+\infty
        \qquad
        \text{ as }\lambda\downarrow\max_{Q_b}\bar a.
    \end{equation}

    \begin{corollary}\label{cor:totally-irrational}
        If $b_1,\ldots,b_n$ are rationally independent, then $H_b=\mathbb{T}^n$,
        $Q_b$ is a single point, and
        \[
        \omega(a,V_b)=\int_{\mathbb{T}^n}a(x)\,\mathrm{d} x.
        \]
        Consequently, $V_b+a+\mathcal{K}$ has a principal eigenvalue
        $\lambda_*=s(V_b+a+\mathcal{K})$ satisfying
        \[
        \lambda_*>\int_{\mathbb{T}^n}a(x)\,\mathrm{d}x.
        \]
    \end{corollary}

    \begin{proof}
        When $b_1,\ldots,b_n$ are rationally independent, the orbit closure is $\mathbb{T}^n$ (see, e.g.,~\cite[Section~1.5]{KatokHasselblatt1995Book}). In this case \(Q_b\) consists of a single point, and the quotient integral is
        \[
        \left(
        \lambda-\int_{\mathbb{T}^n}a(x)\,\mathrm{d} x
        \right)^{-1},
        \]
        which diverges as
        \(\lambda\downarrow\int_{\mathbb{T}^n}a(x)\,\mathrm{d} x\).
        The conclusion therefore follows from
        Theorem~\ref{thm:intro-C}.
    \end{proof}

    Corollary~\ref{cor:totally-irrational} also follows from
    Theorem~\ref{thm:intro-B}. Indeed, the totally irrational translation
    flow on \(\mathbb{T}^n\) is uniquely ergodic with normalized Haar
    measure, and condition \textup{(AC)} holds with density
    \(\rho\equiv1\).

    \begin{corollary}\label{cor:partially-rational}
        Suppose that $\dim H_b=d<n$, and set
        $\bar a_*:=\max_{Q_b}\bar a$.  If, in local coordinates near a
        maximum point $q_0$,
        \[
        \bar a_* - \bar a(q)\le C|q-q_0|^\alpha
        \qquad q\text{ near }q_0
        \]
        for some $C>0$ and $\alpha\ge n-d$, then $V_b+a+\mathcal{K}$ has a
        principal eigenvalue $\lambda_*=s(V_b+a+\mathcal{K})$ satisfying
        \[
        \lambda_*>\bar a_*.
        \]
    \end{corollary}

    \begin{proof}
        Set \(k:=n-d\). Since \(Q_b\) is a \(k\)-dimensional torus and
        \(\eta_b\) is its normalized Haar measure, in local coordinates near
        \(q_0\), the measure \(\eta_b\) is comparable to \(k\)-dimensional
        Lebesgue measure. Moreover,
        \[
        \lambda-\bar a(q)
        \le
        \lambda-\bar a_*+C|q-q_0|^\alpha
        \]
        near \(q_0\). Hence, for some \(c>0\) and \(\delta>0\),
        \[
        \int_{Q_b}\frac{\,\mathrm{d}\eta_b(q)}
        {\lambda-\bar a(q)}
        \ge
        c\int_0^\delta
        \frac{\rho^{k-1}}
        {\lambda-\bar a_*+C\rho^\alpha}\,\mathrm{d}\rho.
        \]
        As \(\lambda\downarrow\bar a_*\), the integral on the right diverges
        whenever \(\alpha\ge k\). Thus
        \eqref{eq:torus-integral-condition} holds, and the conclusion follows
        from Theorem~\ref{thm:intro-C}.
    \end{proof}

    \begin{remark}\label{rem:applications-torus-time-change}
        If $V_b$ is replaced by $fV_b$, where $f\in C^1(\mathbb{T}^n)$ and $f>0$,
        the leaves and the quotient are unchanged.  In the formula for the
        transport growth bound and in the quotient criterion, simply replace
        $\bar a$ by
        \[
        \bar a_f(q)
        =\frac{\int_{H_b}a(x+h)f(x+h)^{-1}\,\mathrm{d} m_{H_b}(h)}
        {\int_{H_b}f(x+h)^{-1}\,\mathrm{d} m_{H_b}(h)}.
        \]
        See Corollary~\ref{cor:time-changed-group}.
    \end{remark}

    We now specialize to $M=\mathbb{T}^2$.  If $b=(1,\alpha)$ with
    $\alpha\notin\mathbb{Q}$, then $H_b=\mathbb{T}^2$, and hence
    \begin{equation}\label{eq:T2-irrational}
        \omega(a,b\cdot\nabla)=\int_{\mathbb{T}^2}a(x,y)\,\mathrm{d} x\,\mathrm{d} y.
    \end{equation}
    Therefore $b\cdot\nabla+a+\mathcal{K}$ has a principal eigenvalue.

    If $b$ has rational direction, an element of $\mathrm{GL}(2,\mathbb{Z})$
    reduces the flow to a nonzero constant multiple of $\partial_x$.
    We use the normalized field $V=\partial_x$.  Its leaves are the
    horizontal circles
    \[
    F_y=\mathbb{T}^1\times\{y\},
    \]
    and the quotient is $Q=\mathbb{T}^1_y$.  Thus
    \begin{equation}\label{eq:T2-rational-average}
        \bar a(y)=\int_0^1a(x,y)\,\mathrm{d} x,
        \qquad
        \omega(a,\partial_x)=\max_{y\in\mathbb{T}^1}\bar a(y).
    \end{equation}
    The operator $\partial_x+a+\mathcal{K}$ has a principal eigenvalue if
    \begin{equation}\label{eq:T2-rational-integral-condition}
        \int_{\mathbb{T}^1}\frac{\,\mathrm{d} y}{\lambda-\bar a(y)}\to+\infty
        \qquad
        \text{ as } \lambda\downarrow\max\bar a.
    \end{equation}
    In particular, a nondegenerate quadratic maximum of $\bar a$ is
    sufficient.

    \subsection{Contrasting effects of transport in a rational direction}

    The following two examples show that transport in a rational direction may
    either create or suppress a principal eigenvalue.

    \subsubsection*{Transport-induced principal eigenvalues}

    The first example shows that, for sufficiently weak nonlocal dispersal, transport in a rational direction can create a principal eigenvalue even though the operator without transport has none.

    Let $M=\mathbb{T}^2$.  Fix $1<\gamma<2$ and choose a nonnegative function
    $p\in C^1(\mathbb{T}^1)$ such that
    \[
    p(x)=|x|^\gamma\quad\text{near }x=0,
    \qquad
    p(x)>0\quad\text{for }x\ne0.
    \]
    Let $0<\sigma<1$ and set
    \begin{equation}\label{eq:example-a}
        a(x,y)=1-p(x)+\sigma\cos(2\pi y).
    \end{equation}
    Then the unique point of maximum is $(0,0)$ and
    \begin{equation}\label{eq:creation-pointwise-asymptotic}
        a_* - a(x,y)\asymp |x|^\gamma+y^2
        \qquad\text{near }(0,0),
    \end{equation}
    where $a_*:=\max_{\mathbb{T}^2}a=1+\sigma$.

    \begin{lemma}\label{lem:finite-no-advection}
        The following estimate holds:
        \begin{equation}\label{eq:finite-zero-advection-integral}
            \int_{\mathbb{T}^2}\frac{\,\mathrm{d} x\,\mathrm{d} y}{a_*-a(x,y)}<+\infty.
        \end{equation}
    \end{lemma}

    \begin{proof}
        Away from $(0,0)$, the denominator is bounded below by a positive
        constant.  Near $(0,0)$, \eqref{eq:creation-pointwise-asymptotic}
        reduces the question to
        \[
        \int_{|x|,|y|<\delta}\frac{\,\mathrm{d} x\,\mathrm{d} y}{|x|^\gamma+y^2}<+\infty.
        \]
        For fixed $x\ne0$,
        \[
        \int_{-\delta}^{\delta}\frac{\,\mathrm{d} y}{|x|^\gamma+y^2}
        \le C |x|^{-\gamma/2}.
        \]
        Since $\gamma<2$, $\int_0^\delta x^{-\gamma/2}\,\mathrm{d} x<\infty$.
        This proves the lemma.
    \end{proof}

    Consider the operator without transport
    \[
    \mathcal{L}_0^\varepsilon u=a(x,y)u+\varepsilon\mathcal{K} u,
    \qquad \varepsilon>0.
    \]
    Now introduce transport in the rational direction
    \(V=\partial_x\) and define the corresponding operator by
    \[
    \mathcal{L}_1^\varepsilon u
    =\partial_xu+a(x,y)u+\varepsilon\mathcal{K}u.
    \]
    The leaf average is
    \begin{equation}\label{eq:creation-leaf-average}
        \bar a(y)=\int_0^1 a(x,y)\,\mathrm{d} x
        =1-\bar p+\sigma\cos(2\pi y),
        \qquad
        \bar p=\int_0^1p(x)\,\mathrm{d} x,
    \end{equation}
    Formula \eqref{eq:creation-leaf-average} gives
    \[
    \bar a_*:=\max_y\bar a(y)=1-\bar p+\sigma,
    \qquad
    \bar a_* - \bar a(y)\asymp y^2
    \quad\text{near }y=0.
    \]
    Therefore
    \begin{equation}\label{eq:creation-quotient-divergence}
        \int_{\mathbb{T}^1}\frac{\,\mathrm{d} y}{\lambda-\bar a(y)}\to+\infty
        \qquad \text{ as }\lambda\downarrow\bar a_*.
    \end{equation}
    The comparison can now be stated as a single theorem.

    \begin{theorem}\label{thm:advection-creates}
        The following statements are valid:
        \begin{enumerate}[label=\textup{(\roman*)},leftmargin=2.6em]
            \item There exists \(\varepsilon_0>0\) such that, for every
            \(0<\varepsilon<\varepsilon_0\), the operator without transport
            $\mathcal{L}_0^\varepsilon$ satisfies
            \begin{equation}\label{eq:creation-threshold-equality}
                s(\mathcal{L}_0^\varepsilon)=a_*=\omega(a,0).
            \end{equation}
            Consequently, $\mathcal{L}_0^\varepsilon$ has no principal
            eigenvalue on \(C(\mathbb{T}^2)\).
            \item For every \(\varepsilon>0\), the operator $\mathcal{L}_1^\varepsilon$ has a principal eigenvalue and
            \(
            s(\mathcal{L}_1^\varepsilon)>\bar a_*.
            \)
        \end{enumerate}
    \end{theorem}

    \begin{proof}
        For \textup{(i)}, let $\mathcal{A}_0$ be the multiplication operator
        $\mathcal{A}_0u=au$.  Then
        \[
        s(\mathcal{A}_0)=a_*=\omega(a,0),
        \]
        and, for $\lambda>a_*$,
        \[
        (R_\lambda^{(0)}u)(x,y)
        :=((\lambda-\mathcal{A}_0)^{-1}u)(x,y)
        =\frac{u(x,y)}{\lambda-a(x,y)}.
        \]
        Set
        \[
        M_0:=\sup_{\lambda>a_*}\sup_{(x,y)\in\mathbb{T}^2}
        \int_{\mathbb{T}^2}\frac{K((x,y),(x',y'))}{\lambda-a(x',y')}\,\mathrm{d} x'\,\mathrm{d} y'.
        \]
        By continuity of \(K\) and
        \eqref{eq:finite-zero-advection-integral},
        \[
        M_0\le \|K\|_\infty
        \int_{\mathbb{T}^2}\frac{\,\mathrm{d} x\,\mathrm{d} y}{a_*-a(x,y)}<\infty.
        \]
        Positivity of the kernel gives
        \[
        \|\mathcal{K}R_\lambda^{(0)}\|\le M_0,
        \qquad \lambda>a_*.
        \]
        Corollary~\ref{cor:critical-coupling} and Remark~\ref{rem:critical-coupling-bounds} yield
        \eqref{eq:creation-threshold-equality} and the absence of a
        principal eigenvalue whenever $\varepsilon M_0<1$.  Thus
        \textup{(i)} holds, for example, with
        $\varepsilon_0=M_0^{-1}$.

        For \textup{(ii)}, the operator is \(\partial_x+a+\varepsilon\mathcal{K}\).  The kernel \(\varepsilon K\) satisfies the same standing kernel hypothesis.  The quotient is \(\mathbb{T}^1_y\), and
        \eqref{eq:creation-quotient-divergence} is exactly the quotient integral
        condition \eqref{eq:T2-rational-integral-condition}.  Theorem~\ref{thm:intro-C} applies.
    \end{proof}

    \subsubsection*{Transport-suppressed principal eigenvalues}

    The second example shows the opposite phenomenon: transport in a rational
    direction can suppress a principal eigenvalue that is present in the purely
    nonlocal problem.

    Fix $0<\beta<1$ and $0<\delta<\tfrac14$.  Choose nonnegative
    functions $p_1,p_2,h_\beta,w\in C(\mathbb{T}^1)$ such that, in the standard
    local coordinate centered at $0$,
    \[
    p_1(x)=|x|^4\quad\text{for }|x|<2\delta,
    \qquad
    p_2(y)=|y|^4\quad\text{for }|y|<2\delta,
    \]
    and
    \[
    h_\beta(y)=|y|^\beta\quad\text{for }|y|<2\delta.
    \]
    Assume that $p_1(x)>0$ for $x\ne0$, while $p_2(y)>0$ and
    $h_\beta(y)>0$ for $y\ne0$.  Finally, choose $w\not\equiv0$ so that
    \[
    w(x)=0\qquad\text{whenever }|x|<2\delta.
    \]
    Since $w\ge0$ and $w\not\equiv0$,
    \[
    \bar w:=\int_0^1 w(x)\,\mathrm{d} x>0.
    \]
    Define
    \begin{equation}\label{eq:suppression-potential}
        a(x,y)=1-p_1(x)-p_2(y)-w(x)h_\beta(y).
    \end{equation}
    Then $a_*:=\max_{\mathbb{T}^2}a=1$, attained only at $(0,0)$.  If
    $|x|,|y|<\delta$, then the support condition gives $w(x)=0$, and
    hence
    \[
    a_* - a(x,y)=|x|^4+|y|^4.
    \]
    Consequently,
    \begin{equation}\label{eq:suppression-pointwise-divergence}
        \int_{\mathbb{T}^2}\frac{\,\mathrm{d} x\,\mathrm{d} y}{a_*-a(x,y)}=+\infty.
    \end{equation}
    On the other hand, the horizontal leaf average satisfies
    \begin{equation}\label{eq:suppression-leaf-asymptotic}
        \begin{aligned}
            \bar a(y)&=\int_0^1 a(x,y)\,\mathrm{d} x
            =1-\bar p_1-p_2(y)-\bar w h_\beta(y),\\
            \bar a_*&:=\max_y\bar a(y)=1-\bar p_1,
            \qquad
            \bar a_* -\bar a(y)\asymp |y|^\beta
            \quad\text{near }0,
        \end{aligned}
    \end{equation}
    where $\bar p_1:=\int_0^1p_1(x)\,\mathrm{d} x$.  By
    \eqref{eq:suppression-leaf-asymptotic} and $\beta<1$,
    \begin{equation}\label{eq:suppression-quotient-integrability}
        \int_{\mathbb{T}^1}\frac{\,\mathrm{d} y}{\bar a_* -\bar a(y)}<+\infty.
    \end{equation}

    \begin{theorem}\label{thm:advection-suppresses}
        The following statements are valid:
        \begin{enumerate}[label=\textup{(\roman*)},leftmargin=2.6em]
            \item For every $\varepsilon>0$, the operator without transport
            \[
            \widetilde{\mathcal{L}}_0^\varepsilon u=a(x,y)u+\varepsilon\mathcal{K} u
            \]
            has a principal eigenvalue on $C(\mathbb{T}^2)$.
            \item There exists $\varepsilon_0>0$ such that, for every
            $0<\varepsilon<\varepsilon_0$, the operator
            \[
            \widetilde{\mathcal{L}}_1^\varepsilon u=\partial_xu+a(x,y)u+\varepsilon\mathcal{K} u
            \]
            satisfies
            \begin{equation}\label{eq:suppression-threshold-equality}
                s(\widetilde{\mathcal{L}}_1^\varepsilon)
                =\bar a_*=\omega(a,\partial_x).
            \end{equation}
            Consequently, $\widetilde{\mathcal{L}}_1^\varepsilon$ has no
            principal eigenvalue on $C(\mathbb{T}^2)$.
        \end{enumerate}
    \end{theorem}

    \begin{proof}
        For \textup{(i)}, the zero-transport resolvent is multiplication by
        $(\lambda-a)^{-1}$.  By \eqref{eq:suppression-pointwise-divergence},
        \[
        \int_{\mathbb{T}^2}\frac{\,\mathrm{d} x\,\mathrm{d} y}{\lambda-a(x,y)}\to+\infty
        \qquad(\lambda\downarrow a_*).
        \]
        Proposition~\ref{prop:divergence-condition} shows that the
        spectral radius is infinite.  Hence Corollary~\ref{cor:critical-coupling}
        gives a principal eigenvalue for every $\varepsilon>0$.

        We prove \textup{(ii)}.  Let $\mathcal{A}_1=\partial_x+a$ and let
        $R_\lambda^{(1)}=(\lambda-\mathcal{A}_1)^{-1}$ for
        $\lambda>\bar a_*$.  We first note that
        \begin{equation}\label{eq:periodic-resolvent-L1-bound}
            \sup_{\lambda>\bar a_*}\int_{\mathbb{T}^2}R_\lambda^{(1)}\mathbf{1}\,\mathrm{d} x\,\mathrm{d} y<+\infty.
        \end{equation}
        Indeed, for fixed $y$, the flow in the $x$ variable has period one and
        average $\bar a(y)$.  The positive resolvent has the representation
        \[
        (R_\lambda^{(1)}\mathbf{1})(x,y)
        =\int_0^\infty
        \exp\left(-\lambda t+\int_0^t a(x+s,y)\,\mathrm{d} s\right)\,\mathrm{d} t,
        \qquad \lambda>\bar a_*,\quad (x,y)\in\mathbb{T}^2.
        \]
        Since $a(\cdot,y)-\bar a(y)$ is one-periodic with zero mean, there is
        a constant $C_0>0$ such that
        \[
        \left|\int_0^t\bigl(a(x+s,y)-\bar a(y)\bigr)\,\mathrm{d} s\right|\le C_0
        \qquad (x,y)\in\mathbb{T}^2,\quad t\ge0.
        \]
        Therefore
        \[
        (R_\lambda^{(1)}\mathbf{1})(x,y)
        \le \frac{e^{C_0}}{\lambda-\bar a(y)},
        \qquad \lambda>\bar a_*.
        \]
        Together with \eqref{eq:suppression-quotient-integrability}, this
        proves \eqref{eq:periodic-resolvent-L1-bound}.

        Set
        \[
        M_1:=\|K\|_\infty
        \sup_{\lambda>\bar a_*}\int_{\mathbb{T}^2}R_\lambda^{(1)}\mathbf{1}\,\mathrm{d} x\,\mathrm{d} y<+\infty.
        \]
        The estimate in Remark~\ref{rem:critical-coupling-bounds} gives
        \[
        \|\mathcal{K}R_\lambda^{(1)}\|\le M_1,
        \qquad \lambda>\bar a_*.
        \]
        For this periodic transport, $\mathcal{A}_1$ is conjugate by a positive
        periodic multiplication operator to $\partial_x+\bar a(y)$.  Functions
        independent of $x$ then show that
        \[
        s(\mathcal{A}_1)=\max_y\bar a(y)=\bar a_*
        =\omega(a,\partial_x).
        \]
        Corollary~\ref{cor:critical-coupling} now yields
        \eqref{eq:suppression-threshold-equality} and the absence of a principal
        eigenvalue whenever $\varepsilon M_1<1$.  Thus \textup{(ii)} holds, for
        example, with $\varepsilon_0=M_1^{-1}$.
    \end{proof}

    \subsubsection*{Discussion}

    These examples show that the criterion is controlled by the singularity on
    the orbit quotient.  Without transport, the criterion depends on $a_*-a$ on
    $\mathbb{T}^2$; under $V=\partial_x$, it depends on $\bar a_*-\bar a$ on the
    one-dimensional quotient $\mathbb{T}^1$.  Leaf averaging can change reciprocal
    integrability in either direction, so transport in a rational direction may
    create or suppress a principal eigenvalue.  In the terminology of
    Corollary~\ref{cor:critical-coupling}, transport in a rational direction
    changes the critical coupling from a positive value to zero in
    Theorem~\ref{thm:advection-creates},
    and from zero to a positive value in
    Theorem~\ref{thm:advection-suppresses}.

    \subsection{Rotations on spheres}

    We first consider axial rotations on $S^2$.  Let $M=S^2\subset\mathbb{R}^3$ and
    write
    \[
    x=\sqrt{1-z^2}\cos\theta,
    \qquad
    y=\sqrt{1-z^2}\sin\theta,
    \qquad
    \theta\in\mathbb{R}/2\pi\mathbb{Z},
    \quad z\in[-1,1].
    \]
    Consider the rotation field
    \[
    V=\partial_\theta.
    \]
    For $|z|<1$, the leaves are the latitude circles
    \[
    F_z=\{(\sqrt{1-z^2}\cos\theta,
    \sqrt{1-z^2}\sin\theta,z):\theta\in\mathbb{R}/2\pi\mathbb{Z}\},
    \]
    while the poles form the two singular leaves $F_1=\{N\}$ and
    $F_{-1}=\{S\}$.  Hence the quotient is
    \[
    Q=[-1,1].
    \]
    The leaf average and transport growth bound are
    \begin{equation}\label{eq:sphere-leaf-threshold}
        \bar a(z)=\frac{1}{2\pi}\int_0^{2\pi}a(\theta,z)\,\mathrm{d}\theta
        \quad(|z|<1),
        \qquad
        \omega(a,\partial_\theta)=\max_{z\in[-1,1]}\bar a(z),
    \end{equation}
    with $\bar a(1)=a(N)$ and $\bar a(-1)=a(S)$.

    \begin{corollary}\label{cor:sphere}
        With $\bar a$ given by \eqref{eq:sphere-leaf-threshold}, set
        $\bar a_*:=\max_{[-1,1]}\bar a$.  If
        \begin{equation}\label{eq:sphere-integral-condition}
            \int_{-1}^{1}\frac{\,\mathrm{d} z}{\lambda-\bar a(z)}\to+\infty
            \qquad
            \text{ as }\lambda\downarrow\bar a_*,
        \end{equation}
        then $\partial_\theta+a+\mathcal{K}$ has a principal eigenvalue
        \[
        \lambda_*=s(\partial_\theta+a+\mathcal{K})>\bar a_*.
        \]
    \end{corollary}

    \begin{proof}
        In the coordinates above, normalized area measure is
        $(4\pi)^{-1}\,\mathrm{d}\theta\,\mathrm{d} z$, so the quotient measure is
        $\frac12\,\mathrm{d} z$.  The conclusion follows from
        Theorem~\ref{thm:intro-C}.
    \end{proof}

    \begin{example}[Interior and polar maxima]
        If $z_0\in(-1,1)$ and
        \[
        \bar a(z)=\bar a_*-c(z-z_0)^2+o((z-z_0)^2),
        \qquad z\to z_0,
        \quad c>0,
        \]
        then \eqref{eq:sphere-integral-condition} holds.  If the north pole
        is the unique maximizer of $\bar a$ and
        \[
        \bar a_* - \bar a(z)\asymp(1-z)^\alpha
        \qquad(z\uparrow1)
        \]
        for some $\alpha>0$, then the integral diverges if and only if
        $\alpha\ge1$.
    \end{example}

    \begin{remark}\label{rem:applications-sphere-time-change}
        If $f\in C^1(S^2)$ and $f>0$, replacing $\partial_\theta$ by
        $f\partial_\theta$ leaves the latitude circles unchanged.  Replace
        $\bar a$ by
        \[
        \bar a_f(z)
        =\frac{\int_0^{2\pi}a(\theta,z)f(\theta,z)^{-1}\,\mathrm{d}\theta}
        {\int_0^{2\pi}f(\theta,z)^{-1}\,\mathrm{d}\theta}.
        \]
        Corollary~\ref{cor:sphere} then holds with $\bar a_f$ (see
        Corollary~\ref{cor:time-changed-group}).
    \end{remark}

    The same quotient picture extends to higher-dimensional spheres.  On
    $S^{2\ell-1}\subset\mathbb{C}^\ell$, let $(\mathbb{R}/2\pi\mathbb{Z})^\ell$ act by
    \[
    (e^{i\theta_1},\ldots,e^{i\theta_\ell})\cdot(z_1,\ldots,z_\ell)
    =(e^{i\theta_1}z_1,\ldots,e^{i\theta_\ell}z_\ell),
    \]
    and consider
    \[
    V_\Omega=\sum_{j=1}^\ell\Omega_j\partial_{\theta_j}.
    \]
    If $\Omega_1,\ldots,\Omega_\ell$ are rationally independent, then the
    quotient is the simplex
    \[
    \Delta_{\ell-1}
    =\left\{\rho\in[0,1]^\ell:\sum_{j=1}^\ell\rho_j=1\right\},
    \qquad \rho_j=|z_j|^2.
    \]
    On the interior stratum the leaf average is
    \[
    \bar a(\rho)
    =\frac{1}{(2\pi)^\ell}
    \int_0^{2\pi}\!\cdots\!\int_0^{2\pi}
    a(\rho,\theta)\,\mathrm{d}\theta_1\cdots\,\mathrm{d}\theta_\ell,
    \qquad
    \rho_j>0,
    \quad \sum_{j=1}^\ell\rho_j=1.
    \]
    On a boundary stratum, one averages only over the active angular
    variables.  The transport growth bound is the maximum of these leaf
    averages on the quotient, and Theorem~\ref{thm:intro-C} gives the
    corresponding quotient criterion.  Thus $S^2$ is an explicit example of
    the same compact-group mechanism.

    \section{Conclusion}

    We have established a necessary and sufficient criterion for the existence of
    a principal eigenvalue for
    \(\mathcal{L}=V+a+\mathcal{K}\) on connected compact Riemannian manifolds.
    The operator \(\mathcal{L}\) has a principal eigenvalue if and only if
    \(s(\mathcal{L})>\omega(a,V)\), equivalently, if
    \(r(R_{\lambda_0}\mathcal{K})>1\) for some
    \(\lambda_0>\omega(a,V)\). Whenever it exists, the principal eigenvalue is
    isolated and algebraically simple, has a strongly positive eigenfunction, and
    is the unique spectral point with maximal real part. Since \(\omega(a,V)\) is
    the maximal invariant average of the potential, the relevant threshold is
    determined by long-time transport dynamics rather than by the pointwise
    maximum of the potential alone.

    The criterion is applied in two dynamical settings. For a uniquely ergodic
    flow, the threshold is the unique invariant average. The 
    resolvent estimate used to prove existence is obtained from either a uniform
    one-sided lower bound on Birkhoff deviations \textup{(BD)} or bounded density
    of the invariant measure with respect to Riemannian volume \textup{(AC)}. For
    flows arising from compact abelian Lie group actions, Haar averaging
    identifies the leaf averages, and divergence of the corresponding scalar
    integral over the orbit quotient gives a sufficient existence condition. The
    quotient formulation permits orbit closures of different dimensions and
    extends to positive time changes after inverse-speed reweighting. In both
    settings, identifying the transport threshold and establishing sufficient
    resolvent growth near it are separate steps.

    The torus constructions show that transport can either promote or suppress the existence of a principal eigenvalue. Let \(\varepsilon_c\) denote the critical coupling strength for \(V+a+\varepsilon\mathcal{K}\), so that a principal eigenvalue exists if and only if \(\varepsilon>\varepsilon_c\). Adding the rational transport \(V=\partial_x\) reduces \(\varepsilon_c\) from a positive value to zero for the first potential, whereas it increases \(\varepsilon_c\) from zero to a positive value for the second. Both transitions stem from how leaf averaging alters reciprocal integrability at the transport threshold, with the suppression case also relying on a periodic resolvent estimate. Thus, the net effect of transport depends on the orbit geometry and the resulting resolvent behavior.
    \section*{Declaration on the Use of AI}
    
    The human authors initiated and led the research program, formulated the mathematical problem, developed the main results, and independently verified all arguments. ChatGPT provided relevant concepts and suggestions and assisted in identifying, summarizing, comparing, and explaining background literature, some of which proved particularly useful. ChatGPT also assisted with mathematical exposition, notation, grammar, and presentation. The final manuscript was carefully revised and approved by the human authors. The authors assume responsibility for all content.       


\begin{thebibliography}{99}

        \bibitem{BandleGonzalezFontelosWolanski2018CPDE}
        C. Bandle, M. del Mar Gonz\'alez, M. A. Fontelos and N. Wolanski,
        A nonlocal diffusion problem on manifolds,
        \emph{Comm. Partial Differential Equations} \textbf{43} (2018), 652--676.

        \bibitem{BaoShen2017PAMS}
        X. Bao and W. Shen,
        Criteria for the existence of principal eigenvalues of time periodic cooperative linear systems with nonlocal dispersal,
        \emph{Proc. Amer. Math. Soc.} \textbf{145} (2017), 2881--2894.

        \bibitem{BatesZhao2007JMAA}
        P. W. Bates and G. Zhao,
        Existence, uniqueness and stability of the stationary solution to a nonlocal evolution equation arising in population dispersal,
        \emph{J. Math. Anal. Appl.} \textbf{332} (2007), 428--440.

        \bibitem{BerestyckiCovilleVo2016JFA}
        H. Berestycki, J. Coville and H.-H. Vo,
        On the definition and the properties of the principal eigenvalue of some nonlocal operators,
        \emph{J. Funct. Anal.} \textbf{271} (2016), 2701--2751.

        \bibitem{BerestyckiHamelNadirashvili2005JEMS}
        H. Berestycki, F. Hamel and N. Nadirashvili,
        The speed of propagation for KPP type problems. I. Periodic framework,
        \emph{J. Eur. Math. Soc.} \textbf{7} (2005), 173--213.

        \bibitem{BerestyckiNirenbergVaradhan1994CPAM}
        H. Berestycki, L. Nirenberg and S. R. S. Varadhan,
        The principal eigenvalue and maximum principle for second-order elliptic operators in general domains,
        \emph{Comm. Pure Appl. Math.} \textbf{47} (1994), 47--92.

        \bibitem{BerestyckiRossi2006JEMS}
        H. Berestycki and L. Rossi,
        On the principal eigenvalue of elliptic operators in \(\mathbb{R}^N\) and applications,
        \emph{J. Eur. Math. Soc.} \textbf{8} (2006), 195--215.

        \bibitem{Bredon1972Book}
        G. E. Bredon,
        \emph{Introduction to Compact Transformation Groups},
        Pure and Applied Mathematics, Vol.~46, Academic Press, 1972.

        \bibitem{Burger1988MZ}
        R. B\"urger,
        Perturbations of positive semigroups and applications to population genetics,
        \emph{Math. Z.} \textbf{197} (1988), 259--272.

        \bibitem{Coville2010JDE}
        J. Coville,
        On a simple criterion for the existence of a principal eigenfunction of some nonlocal operators,
        \emph{J. Differential Equations} \textbf{249} (2010), 2921--2953.

        \bibitem{Coville2013AML}
        J. Coville,
        Singular measure as principal eigenfunction of some nonlocal operators,
        \emph{Appl. Math. Lett.} \textbf{26} (2013), 831--835.

        \bibitem{CovilleDavilaMartinez2013AIHPANL}
        J. Coville, J. D\'avila and S. Mart\'{\i}nez,
        Pulsating fronts for nonlocal dispersion and KPP nonlinearity,
        \emph{Ann. Inst. H. Poincar\'e Anal. Non Lin\'eaire} \textbf{30} (2013), 179--223.

        \bibitem{CovilleHamel2020NA}
        J. Coville and F. Hamel,
        On generalized principal eigenvalues of nonlocal operators with a drift,
        \emph{Nonlinear Anal.} \textbf{193} (2020), 111569.

        \bibitem{Deimling1985Book}
        K. Deimling,
        \emph{Nonlinear Functional Analysis},
        Springer-Verlag, Berlin, Heidelberg, 1985.

        \bibitem{DingLiang2015SIAMJMA}
        W. Ding and X. Liang,
        Principal eigenvalues of generalized convolution operators on the circle and spreading speeds of noncompact evolution systems in periodic media,
        \emph{SIAM J. Math. Anal.} \textbf{47} (2015), 855--896.

        \bibitem{DucrotGrietteLiang2026Preprint}
        A. Ducrot, Q. Griette and X. Liang,
        Existence and multiplicity of principal eigenpairs for nonlocal operator with sign-changing advection in a periodic environment,
        preprint, June 25, 2026.

        \bibitem{EdmundsPotterStuart1972PRSA}
        D. E. Edmunds, A. J. B. Potter and C. A. Stuart,
        Non-compact positive operators,
        \emph{Proc. Roy. Soc. London Ser. A} \textbf{328} (1972), 67--81.

        \bibitem{FengLiRuanXin2024SIAMJMA}
        Y.-X. Feng, W.-T. Li, S. Ruan and M.-Z. Xin,
        Principal spectral theory of time-periodic nonlocal dispersal cooperative systems and applications,
        \emph{SIAM J. Math. Anal.} \textbf{56} (2024), 4040--4083.

        \bibitem{Folland1999Book}
        G. B. Folland,
        \emph{Real Analysis: Modern Techniques and Their Applications},
        2nd ed.,
        Pure and Applied Mathematics,
        John Wiley \& Sons, New York, 1999.

        \bibitem{Folland2016Book}
        G. B. Folland,
        \emph{A Course in Abstract Harmonic Analysis},
        2nd ed., CRC Press, Boca Raton, FL, 2016.

        \bibitem{FreidlinGartner1979SMD}
        M. I. Freidlin and J. G\"artner,
        On the propagation of concentration waves in periodic and random media,
        \emph{Soviet Math. Dokl.} \textbf{20} (1979), 1282--1286.

        \bibitem{Helgason2001Book}
        S. Helgason,
        \emph{Differential Geometry, Lie Groups, and Symmetric Spaces},
        Graduate Studies in Mathematics, Vol.~34,
        American Mathematical Society, Providence, RI, 2001,
        corrected reprint of the 1978 original.

        \bibitem{KangRuan2022MA}
        H. Kang and S. Ruan,
        Principal spectral theory and asynchronous exponential growth for age-structured models with nonlocal diffusion of Neumann type,
        \emph{Math. Ann.} \textbf{384} (2022), 1--49.

        \bibitem{Kato1995Book}
        T. Kato,
        \emph{Perturbation Theory for Linear Operators},
        Classics in Mathematics,
        Springer, Berlin, 1995.

        \bibitem{KatokHasselblatt1995Book}
        A. Katok and B. Hasselblatt,
        \emph{Introduction to the Modern Theory of Dynamical Systems},
        Encyclopedia of Mathematics and its Applications, Vol.~54,
        Cambridge University Press, Cambridge, 1995.

        \bibitem{LiCovilleWang2017DCDS}
        F. Li, J. Coville and X. Wang,
        On eigenvalue problems arising from nonlocal diffusion models,
        \emph{Discrete Contin. Dyn. Syst.} \textbf{37} (2017), 879--903.

        \bibitem{Liang2007CPAM}
        X. Liang and X.-Q. Zhao,
        Asymptotic speeds of spread and traveling waves for monotone semiflows with applications,
        \emph{Comm. Pure Appl. Math.} \textbf{60} (2007), 1--40.

        \bibitem{LiangZhangZhao2019JDE}
        X. Liang, L. Zhang and X.-Q. Zhao,
        The principal eigenvalue for periodic nonlocal dispersal systems with time delay,
        \emph{J. Differential Equations} \textbf{266} (2019), 2100--2124.

        \bibitem{MorroSantAnnaVarandas2020AHP}
        M. Morro, R. Sant'Anna and P. Varandas,
        Ergodic optimization for hyperbolic flows and Lorenz attractors,
        \emph{Ann. Henri Poincar\'e} \textbf{21} (2020), 3253--3283.

        \bibitem{Nussbaum1981LNM}
        R. D. Nussbaum,
        Eigenvectors of nonlinear positive operators and the linear Krein--Rutman theorem,
        in E. Fadell and G. Fournier (eds.),
        \emph{Fixed Point Theory},
        Lecture Notes in Mathematics, Vol.~886,
        Springer, Berlin, 1981, 309--330.

        \bibitem{Pazy1983Book}
        A. Pazy,
        \emph{Semigroups of Linear Operators and Applications to Partial Differential Equations},
        Applied Mathematical Sciences, Vol.~44, Springer, New York, 1983.

        \bibitem{RawalShen2012JDDE}
        N. Rawal and W. Shen,
        Criteria for the existence and lower bounds of principal eigenvalues of time periodic nonlocal dispersal operators and applications,
        \emph{J. Dynam. Differential Equations} \textbf{24} (2012), 927--954.
        
        \bibitem{RodriguezBernalSastreGomez2016PRSE}
        A. Rodr\'iguez-Bernal and S. Sastre-G\'omez,
        Linear non-local diffusion problems in metric measure spaces,
        \emph{Proc. Roy. Soc. Edinburgh Sect. A}
        \textbf{146} (2016), 833--863.
        
        \bibitem{RodriguezBernal2022NA}
        A. Rodr\'iguez-Bernal,
        Principal eigenvalue, maximum principles and linear stability
        for nonlocal diffusion equations in metric measure spaces,
        \emph{Nonlinear Anal.} \textbf{221} (2022), 112887.
        
        
        \bibitem{Rudin1987Book}
        W. Rudin,
        \emph{Real and Complex Analysis},
        3rd ed.,
        McGraw-Hill, New York, 1987.
        
        
        \bibitem{Sadovskaya2009ETDS}
        V. Sadovskaya,
        Dimensional characteristics of invariant measures for circle diffeomorphisms,
        \emph{Ergodic Theory Dynam. Systems} \textbf{29} (2009), 1979--1992.

        \bibitem{Schaefer1974Book}
        H. H. Schaefer,
        \emph{Banach Lattices and Positive Operators},
        Grundlehren der mathematischen Wissenschaften, Vol.~215, Springer, 1974.

        \bibitem{ShenXie2017DCDSB}
        W. Shen and X. Xie,
        Spectral theory for nonlocal dispersal operators with time periodic indefinite weight functions and applications,
        \emph{Discrete Contin. Dyn. Syst. Ser. B} \textbf{22} (2017), 1023--1047.

        \bibitem{ShenZhang2010JDE}
        W. Shen and A. Zhang,
        Spreading speeds for monostable equations with nonlocal dispersal in space periodic habitats,
        \emph{J. Differential Equations} \textbf{249} (2010), 747--795.

        \bibitem{ShenZhang2012PAMS}
        W. Shen and A. Zhang,
        Stationary solutions and spreading speeds of nonlocal monostable equations in space periodic habitats,
        \emph{Proc. Amer. Math. Soc.} \textbf{140} (2012), 1681--1696.

        \bibitem{SuLiLouWang2023JMPA}
        Y.-H. Su, W.-T. Li, Y. Lou and X. Wang,
        Principal spectral theory for nonlocal systems and applications to stem cell regeneration models,
        \emph{J. Math. Pures Appl.} \textbf{176} (2023), 226--281.

        \bibitem{SuWangZhang2023JDE}
        Y.-H. Su, X. Wang and T. Zhang,
        Principal spectral theory and variational characterizations for cooperative systems with nonlocal and coupled diffusion,
        \emph{J. Differential Equations} \textbf{369} (2023), 94--114.

        \bibitem{Thieme1998DCDS}
        H. R. Thieme,
        Remarks on resolvent positive operators and their perturbation,
        \emph{Discrete Contin. Dyn. Syst.} \textbf{4} (1998), 73--90.

        \bibitem{Walters1982Book}
        P. Walters,
        \emph{An Introduction to Ergodic Theory},
        Graduate Texts in Mathematics, Vol.~79, Springer, 1982.

        \bibitem{WangZhang2026JMPA}
        M. Wang and L. Zhang,
        Approximation of the generalized principal eigenvalue of cooperative nonlocal dispersal systems and applications,
        \emph{J. Math. Pures Appl.} \textbf{206} (2026), Paper No. 103835.

        \bibitem{Weinberger2002JMB}
        H. F. Weinberger,
        On spreading speeds and traveling waves for growth and migration models in a periodic habitat,
        \emph{J. Math. Biol.} \textbf{45} (2002), 511--548.

        \bibitem{Zhang2024ANS}
        L. Zhang,
        Principal spectral theory and asymptotic behavior of the spectral bound for partially degenerate nonlocal dispersal systems,
        \emph{Adv. Nonlinear Stud.} \textbf{24} (2024), 1012--1041.

    \end{thebibliography}
\end{document}